\documentclass[letterpaper,10pt,twocolumn]{IEEEtran1}
\IEEEoverridecommandlockouts                                 % This command is only
\renewcommand{\baselinestretch}{1}

\usepackage{amssymb,latexsym,color,amsmath,pifont,epsfig,graphicx}
\usepackage{setspace,cite}
\usepackage{algorithm}
\usepackage{algorithmic}
\usepackage{soul}
\usepackage{graphics}
\usepackage{multirow}
\usepackage{bm}
\usepackage{rotating}
\usepackage{epsfig}
\usepackage{amssymb,latexsym,color,amsmath,pifont,epsfig,mathtools}
\usepackage{graphicx,ctable,booktabs}
\usepackage{indentfirst}
\usepackage{dsfont}
\usepackage{url}

\DeclareMathOperator*{\minimize}{minimize}

\usepackage[subrefformat=parens,labelformat=parens]{subfig}

\newcommand{\DefinedAs}[0]{\mathrel{\mathop:}=}

\newcommand{\bbR}{\mathbb{R}}
\newcommand{\tc}{\textcolor}

\newcommand{\one}{\mathds{1}}

\newcommand{\card}{\mathbf{card}}
\newcommand{\mrd}{\mathrm{d}}
\newcommand{\mre}{\mathrm{e}}

\usepackage{color}

\newcommand{\enma}[1]   {\ensuremath{#1}}

\newcommand{\non}{\nonumber}

\newcommand{\beq}{\begin{equation}}
\newcommand{\eeq}{\end{equation}}
\newcommand{\bseq}{\begin{subequations}}
\newcommand{\eseq}{\end{subequations}}
\newcommand{\beqn}{\begin{eqnarray}}
\newcommand{\eeqn}{\end{eqnarray}}
\newcommand{\ba}{\begin{array}}
\newcommand{\ea}{\end{array}}
\newcommand{\bct}{\begin{center}}
\newcommand{\ect}{\end{center}}
\newcommand{\btab}{\begin{tabular}}
\newcommand{\etab}{\end{tabular}}
\newcommand{\btmz}{\begin{itemize}}
\newcommand{\etmz}{\end{itemize}}
\newcommand{\benum}{\begin{enumerate}}
\newcommand{\eenum}{\end{enumerate}}

\newcommand{\trace}     {\enma{\mathrm{trace}}}

\newcommand{\matbegin}{
        \left[
}
\newcommand{\matend}{
        \right]
}

\newcommand{\tbo}[2]{
  \matbegin \begin{array}{c}
       #1 \\ #2
       \end{array} \matend }

\newcommand{\obt}[2]{
  \matbegin \begin{array}{cc}
       #1 & #2
       \end{array} \matend }
\newcommand{\obth}[3]{
  \matbegin \begin{array}{ccc}
       #1 & #2 & #3
       \end{array} \matend }
\newcommand{\tbt}[4]{
  \matbegin \begin{array}{cc}
       #1 & #2 \\ #3 & #4
       \end{array} \matend }

\newcommand{\be}{\begin{equation}}
\newcommand{\ee}{\end{equation}}

\newcommand{\cplxs}{ C\kern -.35em \rule{0.03 em}{.7 ex}~   }

\def\complex{\hbox{C\kern -.45em \rule{0.03 em}{1.5 ex}}~}

\newcommand{\bi}{\begin{itemize}}
\newcommand{\ei}{\end{itemize}}

\begin{document}

\title{Input-output analysis and decentralized optimal control of inter-area oscillations in power systems}

\author{Xiaofan Wu, Florian D\"orfler, and Mihailo R.\ Jovanovi\'c% <-this % stops a space
\thanks{Financial support from ETH Z\"urich startup grants, the University of Minnesota Informatics Institute Transdisciplinary Faculty Fellowship, and the National Science Foundation under award ECCS-1407958 is gratefully acknowledged.}
\thanks{Xiaofan Wu and Mihailo R.\ Jovanovi\'c are with the Department of Electrical and Computer Engineering, University of Minnesota, Minneapolis, MN 55455.  Emails: {\tt [wuxxx836,mihailo]@umn.edu}. Florian D\"{o}rfler is with the Automatic Control Laboratory at ETH Z\"urich, Switzerland. Email: \texttt{dorfler@control.ee.ethz.ch}.}}

\maketitle

\begin{abstract}
Local and inter-area oscillations in bulk power systems are typically identified using spatial profiles of poorly damped modes, and they are mitigated via carefully tuned decentralized controllers. In this paper, we employ non-modal tools to analyze and control inter-area oscillations. Our input-output analysis examines power spectral density and variance amplification of stochastically forced systems and offers new insights relative to modal approaches. To improve upon the limitations of conventional wide-area control strategies, we also study the problem of signal selection and optimal design of sparse and block-sparse wide-area controllers. In our design, we preserve rotational symmetry of the power system by allowing only relative angle measurements in the distributed controllers. For the IEEE 39 New England model, we examine performance tradeoffs and robustness of different control architectures and show that optimal retuning of fully-decentralized control strategies can effectively guard against local and inter-area oscillations.
\end{abstract}

	\vspace*{-2ex}
\section{Introduction}
\label{Section: Introduction}
% -----------------------------------------------------------------------------------------------------%

Inter-area oscillations in bulk power systems are associated with the dynamics of power transfers and involve groups of synchronous machines that oscillate relative to each other. These system-wide oscillations arise from modular network topologies, heterogeneous machine dynamics, adversely interacting controllers, and large inter-area power transfers. With increased system loads and deployment of renewables in remote areas, long-distance power transfers will eventually outpace the addition of new transmission facilities. This induces severe stress and performance limitations on the transmission network and may even cause instability and outages~\cite{FA-CDM-ID-PS-GS-HE-ZJ:99}.

Traditional analysis and control of inter-area oscillations is based on modal approaches~\cite{KP-NM-TD:10,RL:98}. Typically, inter-area oscillations are identified from the spatial profiles of eigenvectors and participation factors of poorly damped modes~\cite{GR:96,MK-GJR-PK:91}, and they are damped via decentralized controllers, whose gains are carefully tuned using root locus~\cite{MK-GJR-SM-PK:92,NM-LTGL:89}, pole placement~\cite{othman1989design}, adaptive~\cite{malik2014experimental}, robust~\cite{zhu2003robust}, and optimal~\cite{JAT-LS:14} control strategies. To improve upon the limitations of decentralized control, recent research centers at distributed wide-area control strategies that involve the communication of remote signals~\cite{JX-FW-CYC-KPW:06,KS-SNS-SCS:06}. The wide-area control signals are typically chosen to maximize modal observability metrics \cite{AH-KI:02,LPK-RS-BCP:12}, and the control design methods range from root locus criteria to robust and optimal control approaches~\cite{GEB-WS-JHC-TGN-NM:00,YZ-AB:08,MZ-LM-PK-CR-GA:05}.

The spatial profiles of the inter-area modes together with modal controllability and observability metrics were previously used to indicate which wide-area links need to be added and how supplemental damping controllers have to be tuned. Here, we depart from the conventional modal approach and propose a novel methodology for analysis and control of inter-area oscillations. In particular, we use input-output analysis to study oscillations in stochastically forced power systems. A similar approach was recently employed to quantify performance of consensus and synchronization networks~\cite{bamjovmitpat12,farlinjovTAC14sync}.

To identify wide-area control architectures and design optimal sparse controllers, we invoke the paradigm of sparsity-promoting optimal control~\cite{farlinjovACC11,linfarjovACC12,linfarjovTAC13admm,wujovACC14}. Recently, this framework was successfully employed for wide-area control of power systems~\cite{dorjovchebulACC13,dorjovchebulTPS14,SS-UM-FA:14,SP-JAT-RI:15}. Here, we follow the formulation developed in~\cite{wujovACC14} and find a linear state feedback that {simultaneously optimizes a quadratic optimal control criterion (associated with incoherent and poorly damped oscillations) and {induces}} a sparse control architecture. The main novel contributions of our control design approach are highlighted below. We improve the previous results~\cite{dorjovchebulACC13,dorjovchebulTPS14,SS-UM-FA:14,SP-JAT-RI:15} at two levels: first, we preserve rotational symmetry of the original power system by allowing only relative angle measurements in the distributed controller, and, second, we allow identification of block-sparse control architectures, where local information associated with a subsystem is either entirely used (or discarded) for control.

We illustrate the utility of our approach using the IEEE 39 New England model~\cite{JHC-KWC:92}. We show how different sparsity-promoting penalty functions can be used to achieve a desired balance between closed-loop performance and communication complexity. In particular, we demonstrate that the addition of certain long-range communication links and careful retuning of the local controllers represent an effective means for improving system performance. For the New England model, it turns out that properly retuned and {\em fully-decentralized\/}  controllers can perform almost as well as the optimal centralized controllers. Our results thus provide a constructive answer to the much-debated question of whether locally observable oscillations in a power network are also locally controllable~\cite{Eliasson1992damping}.

The remainder of the paper is organized as follows. In Section~\ref{Section: Problem formulation}, we briefly summarize the model, highlight causes of inter-area oscillations, and provide background on input-output analysis of power systems. In Section~\ref{Section:control}, we formulate sparse and block-sparse optimal wide-area control problems under the relative angle measurement restriction. In Section~\ref{Section:Example}, we apply our sparse controllers to the IEEE $39$ New England power grid and compare performance of open- and closed-loop systems. Finally, in Section~\ref{Section:Conclusion}, we conclude the paper.

\section{Background on power system oscillations}
\label{Section: Problem formulation}

\subsection{Modeling and control preliminaries}
\label{Section:Motivation}

A power network is described by a nonlinear system of differential-algebraic equations. Differential equations govern the dynamics of generators and their controllers, and the algebraic equations describe quasi-stationary load flow and circuitry of generators and power electronics~\cite{kundur1994power}. A linearization around a stationary operating point and elimination of the algebraic equations yield a linearized state-space model
	\be
	\dot x
    \; = \;
    {A} \, x
    \; + \;
    {B}_1 \, d
    \; + \;
    {B}_2 \, u.
	\label{eq.ss-ol}	
%	\non
	\ee
Here, $x(t) \in \bbR^n$ is the state, $u(t) \in \bbR^m$ is the generator excitation control input, and $d(t) \in \bbR^p$ is the stochastic disturbance which may arise from power imbalance and uncertain load demands~\cite{kundur1994power}. For example, the choice $B_1= B_2$ can be used to quantify and mitigate the impact of noisy or lossy communication among spatially distributed controllers~\cite{dorjovchebulTPS14}.

The dominant electro-mechanical dynamics of a power system are given by the linearized {\em swing equations}~\cite{kundur1994power},
	\be
%		M \, \ddot{\theta} \; + \;  D \, \dot{\theta} \; + \; L \, \theta
	M_{i} \, \ddot{\theta}_{i} \; + \;  D_{i} \, \dot{\theta}_{i} \; + \; \sum\nolimits_{j} L_{ij} \, (\theta_{i} \, - \, \theta_{j})
	\;=\;
	0.
	% \label{eq.swing}
	\non
	\ee
These equations are obtained by neglecting fast electrical dynamics and eliminating the algebraic load flow. Here, $\theta_{i}$ and $\dot{\theta}_{i}$ are the rotor angle and frequency of generator $i$, $M_{i}$ and $D_{i}$ are the generator inertia and damping coefficients, and $L_{ij}$ is the $(i,j)$ element of the network susceptance matrix indicating the interactions between generators $i$ and $j$~\cite{dorjovchebulTPS14}. Even though the swing equations do not fully capture complexity of power systems, they nicely illustrate the causes of inter-area oscillations: Inter-area oscillations originate from sparse links between densely connected groups of generators (so-called areas). These areas can be aggregated into coherent groups of machines which swing relative to each other using the slow coherency theory~\cite{JHC-PK:85,DR-FD-FB:12q}. Our goal is to design wide-area controllers to suppress inter-area oscillations.

Under a linear state-feedback,
	\[
	u ~=~ -Kx
	\]
the closed-loop system takes the form
    \be
    \ba{rcl}
    \dot{x}
    & \!\! = \!\! &
    ( A \, - \, B_2 K ) \, x
    \; + \;
    B_1 \, d \\[0.15cm]
    z
    & \!\! = \!\! &
    \tbo{z_1}{z_2}
    \; = \;
    \tbo{{Q}^{1/2}}{-{R}^{1/2} \, K} x
    \ea
    \label{eq.ss}
    \ee
where $z$ is a performance output with state and control weights ${Q}$ and ${R}$. We choose $R$ to be the identity matrix and a state objective that quantifies a desired potential energy and the kinetic energy stored in the electro-mechanical dynamics,
    \[
    x^T Q \, x
    \; = \;
    \theta^T Q_{\theta} \, \theta
    \; + \;
    \dfrac{1}{2} \, {\dot{\theta}}^T M \, \dot{\theta}.
    \]
Here, $M = \mathrm{diag} \, (M_i)$ is the inertia matrix and the matrix $Q_{\theta}$ penalizes the deviation of angles from their average
	$
	\bar{\theta} (t)
	\DefinedAs
	(1/N) \, {\one}^T  \theta (t),
	$
    \be
    \ba{rcl}
    Q_{\theta}
    & \!\! = \!\! &
    I \,-\, (1/N) \, \mathds{1} {\mathds{1}}^T
    \ea
    \label{eq.Lunif}
    \ee%
where $N$ is the number of generators and $\mathds{1}$ is the vector of all ones. In a power system without a slack bus, the generator rotor angles are only defined in a relative frame of reference, as can be observed in the swing equations. Thus, they can be rotated by a uniform amount without changing the fundamental dynamics~\eqref{eq.ss-ol}. We preserve this rotational symmetry and study problems in which only differences between the components of the vector $\theta (t) \in \bbR^N$ enter into~\eqref{eq.ss}.
As a result of the rotational symmetry, both the open-loop $A$-matrix and the performance weight $Q_{\theta}$ have an eigenvalue at zero which characterizes the mean of all rotor angles.

By expressing the state vector as
	\be
	x(t)
	\; \DefinedAs \,
	\tbo{\theta(t)}{r(t)}
	\, \in \,
	\bbR^n
	\non
	%\label{eq.svec}
	\ee
where $r (t) \in \bbR^{n-N}$ represents the rotor frequencies and additional states that account for fast electrical dynamics, we arrive at the structural constraints on the matrices in~\eqref{eq.ss},
	\be
	A \tbo{\mathds{1}}{0}
	\; = \;
	0,
	~~
	Q_{\theta} \, \mathds{1}
	\; = \;
	0,
	~~
	K \tbo{\mathds{1}}{0}
	\; = \;
	0.
	\label{eq.assume}
	\non
	\ee
	
In earlier work~\cite{dorjovchebulACC13,dorjovchebulTPS14}, we have removed this rotational symmetry by adding a small {\em regularization\/} term to the diagonal elements of the matrix $Q_{\theta}$. This has resulted in a controller that requires the use of absolute angle measurements to stabilize the average rotor angle. Such a regularization induces a slack bus (a reference generator with a fixed angle) and thereby alters the structure of the original power system.

In this paper, we preserve the natural rotational symmetry by restricting our attention to relative angle measurements. This requirement implies that the average rotor angle has to remain invariant under the state feedback $u=-Kx$. To cope with these additional structural constraints, the sparsity-promoting approach of~\cite{linfarjovTAC13admm} has been recently augmented in~\cite{wujovACC14}.

To eliminate the average-mode $\bar{\theta}$ from~\eqref{eq.ss} we introduce the following coordinate transformation~\cite{wujovACC14},
	\be
	x
	\; = \,
	\tbo{\theta}{r}
	\, = \,
	\underbrace{\tbt{U}{0}{0}{I}}_{T}
	\xi
	~ + \;
	\tbo{\one}{0}
	\bar{\theta}
	\label{eq.x-xi}
	\ee
where the columns of the matrix $U \in \mathbb{R}^{N \times (N-1)}$ form an orthonormal basis that is orthogonal to $\textup{span} \, (\one)$. For example, these columns can be obtained from the $(N-1)$ eigenvectors of the matrix $Q_{\theta}$ in~\eqref{eq.Lunif} that correspond to the non-zero eigenvalues. In the new set of coordinates,
	$
	\xi(t) = T^T x(t) \in \bbR^{n-1},
	$
the closed-loop system takes the form
    \be
    \ba{rcl}
    \dot{\xi}
    & \!\! = \!\! &
    ( \bar{A} \, - \, \bar{B}_2 F ) \, \xi
    \; + \;
    \bar{B}_1 \, d \\[0.15cm]
    z
    & \!\! = \!\! &
    \tbo{z_1}{z_2}
    \; = \;
    \tbo{\bar{Q}^{1/2}}{-R^{1/2} \, F} \xi
    \ea
    \label{eq.ss2cl}
    \ee
where
    \be
  	\bar{A}
	\; \DefinedAs \;
	T^T {A} \, T,
	~~
	\bar B_i
	\; \DefinedAs \;
	T^T {B}_i,
	~~
	\bar Q^{1/2}
    	\; \DefinedAs \;
    	{Q}^{1/2} \, T.
    \label{eq.ss_closed}
    \non
    \ee
The feedback matrices $K$ and $F$ (in the original $x$ and new $\xi$ coordinates, respectively) are related by
	\[
	F
	\; = \;
	K \, T
	~~
	\Leftrightarrow
	~~
	K
	\; = \;
	F \, T^T.
	\]

Because of a marginally stable average mode, the matrix $A$ in~\eqref{eq.ss} is not Hurwitz. The coordinate transformation~\eqref{eq.x-xi} eliminates the average angle $\bar{\theta}$ from~\eqref{eq.ss}, thereby leading to~\eqref{eq.ss2cl} with Hurwitz $\bar{A}$. In the presence of stochastic disturbances, $\bar{\theta} (t)$ drifts in a random walk. Since $\bar{\theta}$ is not observable from the performance output $z$ (which quantifies the mean-square deviation from angle average, kinetic energy, and control effort), $z$ has a finite steady-state variance. This variance is determined by the square of the $\mathcal{H}_2$ norm of system~\eqref{eq.ss2cl}.

	\vspace*{-2ex}
\subsection{Power spectral density and variance amplification}
\label{Section:PSD}

The conventional analysis of inter-area oscillations in power systems is based on spatial profiles of eigenvectors and participation factors of poorly damped modes. Similarly, traditional control design builds on a modal perspective~\cite{GR:96,MK-GJR-PK:91}. In systems with non-normal $A$-matrices, modal analysis may lead to misleading conclusions about transient responses, amplification of disturbances, and robustness margins~\cite{treemb05,jovbamJFM05,liejovkumJFM13}. Non-normal matrices are common in power systems; such matrices do not have orthogonal eigenvectors and they cannot be diagonalized via unitary coordinate transformations.

In what follows, we utilize an approach that offers additional and complementary insights to modal analysis.  This approach is based on the input-output analysis, where the input $d$ is the source of excitation and the output $z$ is the quantity that we care about. In stochastically forced systems, input-output analysis amounts to the study of power spectral density and variance amplification. Our approach builds on the $\mathcal{H}_2$ paradigm~\cite{dulpag00}, which analyzes and mitigates amplification of white stochastic disturbances.

We next provide a brief overview of the power spectral density and variance amplification analyses of linear dynamical systems. Let $H(j \omega)$ denote the frequency response of~\eqref{eq.ss2cl},
	\[
	z (j \omega)
	\; = \;
	H (j \omega) \, d (j \omega).
	\]
The Hilbert-Schmidt norm determines the power spectral density of $H(j \omega)$,
    \be
    {\| H (j \omega) \|}^2_{\rm HS}
    \, = \,
    {\rm{trace}} \left( H (j \omega) \, H^* (j \omega) \right) \\[0.15cm]
    \, = \,
    \displaystyle{\sum} \, \sigma_i^2 (H (j \omega) )
    %\label{eq.HS}
    \non
    \ee
where $\sigma_i$'s are the singular values of the matrix $H(j \omega)$. The $\mathcal{H}_2$ norm quantifies the steady-state variance (energy) of the output $z$ of stochastically forced system~\eqref{eq.ss2cl}. It is obtained by integrating the power spectral density over all frequencies~\cite{dulpag00},
    \be
    {\| H \|}^2_2
    \; \DefinedAs \,
    \displaystyle{\lim_{t \, \to \, \infty}}
    \mathbf{E}
    \left(
    z^T (t) \, z (t)
    \right)
    \, = \;
    \dfrac{1}{2 \pi}
    \, \displaystyle{\int_{- \infty}^{\infty}}
    {\| H (j \omega) \|}^2_{\rm HS} \, \mrd \omega
    \non
    \ee
where $\mathbf{E}$ is the expectation operator. Equivalently, the matrix solution $X$ to the Lyapunov equation,
    \be
    ( \bar{A} \, - \, \bar{B}_2 F ) \, X
    \; + \;
    X \, ( \bar{A} \, - \, \bar{B}_2 F)^T
    \; = \;
    - \bar{B}_1 \bar{B}^T_1
    % \label{eq.lyap}
    \non
    \ee
can be used to compute the $\mathcal{H}_2$ norm~\cite{dulpag00},
     \be
     \ba{rcl}
     J (F)
     \; \DefinedAs \;
    \| H \|^2_2
    & \!\! = \!\! &
    \trace \left( X \, ( {\bar{Q}} \;+\; F^T {R} \, F ) \right)
    \\[0.15cm]
    & \!\! = \!\! &
    \trace \left( Z_1 \right)
    \, + \;
    \trace \left( Z_2 \right).
    \ea
    \label{eq.H2}
    \ee
Here, $X$ is the steady-state covariance matrix of the state $\xi$ in~\eqref{eq.ss2cl},
	$
	X
	\DefinedAs
	\lim_{t \, \to \, \infty}
	 \mathbf{E}
	 \,
    	(
    	\xi (t) \, \xi^T (t)
    	),
	$
and the covariance matrices of the outputs $z_1$ and $z_2$ are determined by
	\[
	\ba{rrl}
	Z_1
	& \!\! \DefinedAs \!\! &
	\displaystyle{\lim_{t \, \to \, \infty}}
    	\mathbf{E}
    	\left(
    	z_1 (t) \, z_1^T (t)
    	\right)
	\; = \;
	{\bar{Q}}^{1/2} X \, {\bar{Q}}^{1/2}
	\\[0.15cm]
	Z_2
	& \!\! \DefinedAs \!\! &
	\displaystyle{\lim_{t \, \to \, \infty}}
    	\mathbf{E}
    	\left(
    	z_2 (t) \, z_2^T (t)
    	\right)
	\; = \;
	R^{1/2} F \, X \, F^T R^{1/2}.
	\ea
	\]
Note that $\trace \, (Z_1)$ and $\trace \, (Z_2)$ quantify the system's kinetic and potential energy and the control effort, respectively. In particular, the eigenvalue decomposition of the matrix $Z_1$,
	\be
    Z_1
    \; = \;
    \displaystyle{\sum} \; \lambda_i \, y_i \, y^T_i
     % \label{eq.eig}
     \non
    \ee
determines contribution of different {\em orthogonal\/} modes $y_i$ to the kinetic and potential energy in statistical steady-state. The total energy is given by $\trace \, (Z_1)$, i.e., the sum of the eigenvalues $\lambda_i$ of the covariance matrix $Z_1$. Each mode $y_i$ contributes $\lambda_i$ to the variance amplification and the spatial structure of the most energetic mode is determined by the principal eigenvector $y_1$ of the matrix $Z_1$.

	\vspace*{-1ex}
\section{Sparse and block-sparse optimal control}
\label{Section:control}

In this section, we study the problem of optimal signal selection and optimal design of wide-area controllers. We approach this problem by invoking sparsity-promoting versions of the standard ${\cal H}_2$ optimal control formulation. We build on the framework developed in~\cite{farlinjovACC11,linfarjovACC12,linfarjovTAC13admm,wujovACC14} which is aimed at finding a state feedback that simultaneously optimizes the closed-loop variance and induces a sparse control architecture. This is accomplished by introducing additional regularization terms to the optimal control problem. These serve as proxies for penalizing the number of communication links in the wide-area controller, thereby inducing a sparse control architecture.

	\vspace*{-2ex}
\subsection{Elementwise sparsity}
\label{Section:sparse}

As shown in Section~\ref{Section:PSD}, the $\mathcal{H}_2$ norm of system~\eqref{eq.ss2cl} is determined by~\eqref{eq.H2}. While the $\mathcal H_{2}$ performance is expressed in terms of the feedback matrix $F$ in the new set of coordinates, it is necessary to enhance sparsity of the feedback matrix $K$ in the physical domain. A desired tradeoff between the system's performance and the sparsity of $K$ is achieved by solving the regularized optimal control problem~\cite{wujovACC14},
    \be
    \ba{rll}
    &
    \minimize\limits_{F, \, K}
    &
    \quad J(F)
    \; + \;
    \gamma \, g(K) \\[0.15cm]
    &
    \rm{subject~~to}
    &
    \quad F \, T^T
    \; - \;
    K
    \; = \;
    0.
    \ea
    \label{eq.ADMM1}
    \ee
The regularization term in~\eqref{eq.ADMM1} is given by the weighted $\ell_1$-norm of $K$,
    \be
    g(K)
    \; \DefinedAs \;
    \displaystyle{\sum_{i, \, j}} \; W_{ij} \, | K_{ij} |
    \non
    \ee
which is an effective proxy for inducing elementwise sparsity~\cite{canwakboy08}. The weights $W_{ij}$'s are updated iteratively using the solution to~\eqref{eq.ADMM1} from the previous iteration; see~\cite{canwakboy08} for details. In~\eqref{eq.ADMM1}, $\gamma$ is a fixed positive scalar that characterizes the emphasis on the sparsity level of the feedback matrix $K$. A larger value of $\gamma$ introduces a sparser feedback gain $K$ at the expense of degrading the closed-loop performance.

We solve the optimal control problem~\eqref{eq.ADMM1} for different values of the positive regularization parameter $\gamma$ via the alternating direction method of multipliers; see~\cite{linfarjovTAC13admm,wujovACC14} for algorithmic details. This allows us to identify a parameterized family of distributed control architectures that strikes an optimal balance between competing performance and sparsity requirements.

	\vspace*{-2ex}
\subsection{Block sparsity}
\label{Section:block}

In power systems, only rotor angle differences enter into the dynamics and information about absolute angles is not available. It is thus advantageous to treat rotor angles separately from the remaining states in the control design. We partition $K$ conformably with the partition of the state vector $x$,
    \[
    K \; = \, \obt{K_\theta}{K_r}
    \]
where $K_\theta$ and $K_r$ are the feedback gains acting on the rotor angles and the remaining states, respectively.

The actuators in wide-area control range from Power System Stabilizers (PSSs) to power electronics devices (FACTS) to HVDC links. While our design methodology is general, in the sequel we restrict our presentation to PSSs. For PSSs the control action is usually formed in a fully-decentralized fashion using local measurements of frequencies and power injections. We represent the vector $r$ as
	\[
	r
	\; = \,
	\obth{r_1^T}{\cdots}{r_N^T}^T
	\]%
where $r_i$ is the vector of states of the controlled generator $i$ (modulo angles). If $K_r$ is partitioned conformably with the partition of the vector $r$, then the block-diagonal elements of $K_r$ provide a means for retuning the local control action. Since $r_i$ is readily available to the controller of generator $i$, in what follows we do not introduce sparsity-promoting penalty on the block-diagonal elements of $K_r$. On the other hand, there are many options for treating the components of $K_r$ associated with the states of other generators. We next illustrate three possible options.

\begin{figure}[!htb]
    \centering
    \subfloat[elementwise]
	    {
    \includegraphics[width=0.30\textwidth]{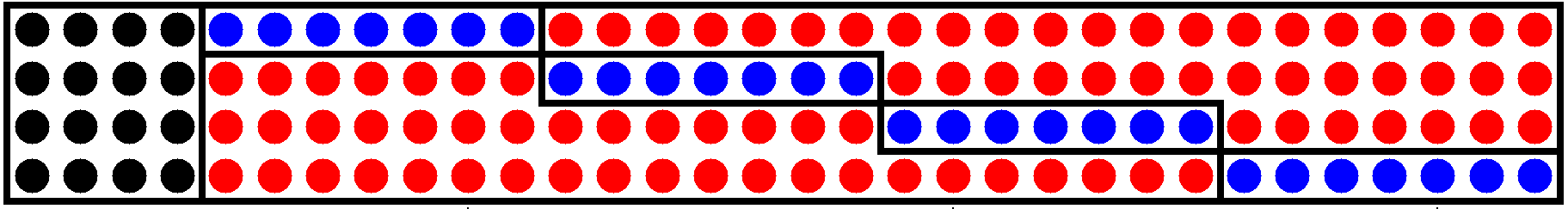}
    \label{fig.block_element}
	    }
	    \\
    \subfloat[group states of individual generators]
	    {
    \includegraphics[width=0.30\textwidth]{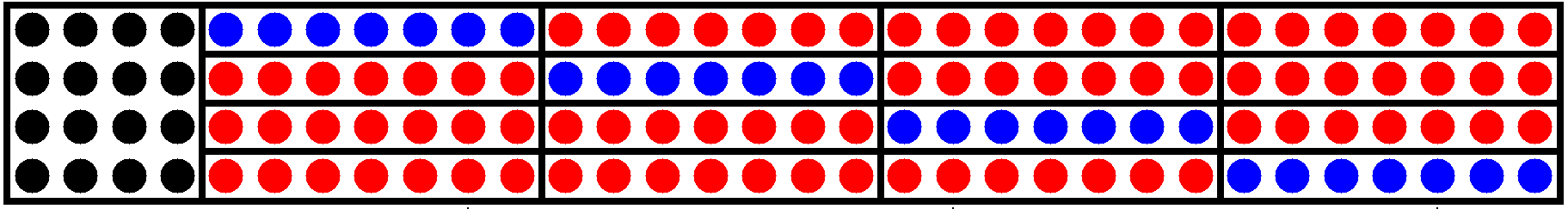}
    \label{fig.block_generator}
	    }
	    \\
    \subfloat[group states of all other generators]
	    {
    \includegraphics[width=0.30\textwidth]{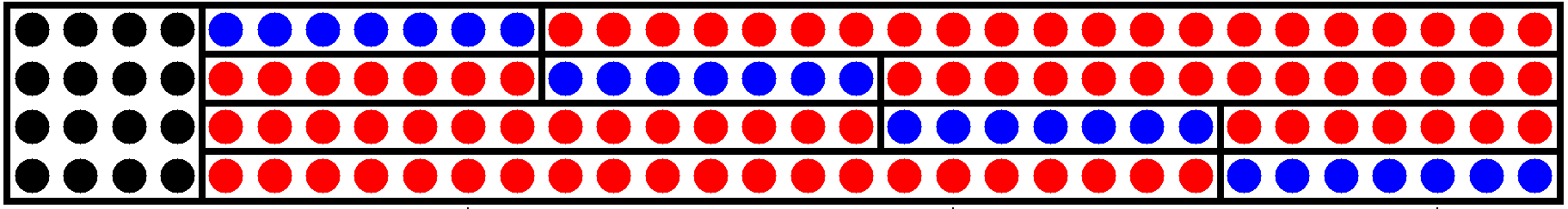}
    \label{fig.block_row}
    }
    \caption{Block structure of the feedback matrix $K$. \tc{black}{\large $\bullet$} denote relative angle feedback gains, \tc{blue}{\large $\bullet$} and \tc{red}{\large $\bullet$} represent local and inter-generator frequency and PSS gains, respectively.}
    \label{fig.block}
    \end{figure}

Consider a system of four generators with controllers. The states of each controlled generator are given by angle, frequency, fluxes, and excitation control system; see Fig.~\ref{fig.block}. Sparsity of the inter-generator control gains can be enhanced either via elementwise or group penalties. Inter-generator information exchange can be treated with an elementwise penalty in the same way as in Section~\ref{Section:sparse}; see Fig.~\ref{fig.block_element} for an illustration. On the other hand, group penalties~\cite{yuan2006model} can be imposed either on the states of individual generators or on the states of all other generators; cf.\ Figs.~\ref{fig.block_generator} and~\ref{fig.block_row}.

\begin{figure}[b]
    \centering
    \begin{tabular}{c}
    \includegraphics[width=0.25\textwidth]{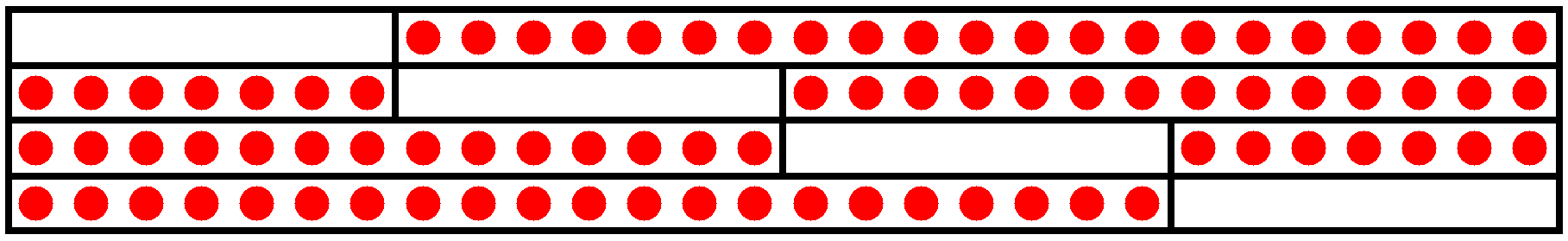}
    \end{tabular}
    \caption{Structural identity matrix $I_s$ with \tc{red}{\large $\bullet$} representing locations of $1$'s.}
    \label{fig.IS}
    \end{figure}

The above objectives can be accomplished by solving the sparsity-promoting optimal control problem
    \be
    \!\!\!\!
    \ba{rll}
    &
    \minimize
    &
    J(F)
    \; + \;
    \gamma_\theta \, g_\theta (K_\theta)
    \; + \;
    \gamma_r \, g_r (K_r)
    \\[0.25cm]
    &
    \rm{subject~to}
    &
    F \, T^T
    \; - \;
    \obt{K_\theta}{K_r}
    \; = \;
    0
    \ea
    \label{eq.ADMM2}
    \ee
where
    \begin{subequations}
	\label{eq.gr}
	\be
    g_\theta (K_\theta)
    \; \DefinedAs \;
    \displaystyle{\sum_{i, \, j}} \; W_{ij} \, | {K_\theta}_{ij} |.
    \label{eq.gtheta}
    \ee%
On the other hand, for the three cases discussed and illustrated in Fig.~\ref{fig.block} the corresponding regularization functions are
	\begin{eqnarray}
	\!\!\!\!\!\!\!\!\!
	g_{r1} (K_r)
	& \!\!\! \DefinedAs \!\!\! &
    	\displaystyle{\sum_{i, \, j}}
	\; W_{ij} \,
	| \left( I_s \, \circ \, K_r \right)_{ij} |
	\label{eq.gr1}
	\\[0.1cm]
	\!\!\!\!\!\!\!\!\!
	g_{r2} (K_r)
	& \!\!\! \DefinedAs \!\!\! &
   	\displaystyle{\sum_{i \, \neq \, k}}
	\, \beta_{ik} \, W_{i k} \,
	|| \, \mre_i^T \! \left( I_s \, \circ \, K_r \right) \circ v_k^T ||_2
    	\label{eq.gr2}
   	 \\[0.1cm]
	 \!\!\!\!\!\!\!\!\!
	 g_{r3} (K_r)
	& \!\!\! \DefinedAs \!\!\! &
    	\displaystyle{\sum_i}
	\, \beta_{i} \, W_i \,
	|| \, \mre_i^T \! \left( I_s \, \circ \, K_r \right) ||_2
    \label{eq.gr3}
	\end{eqnarray}
where
	$
	 i = \{1,\cdots,m\},
   	$
	$
	j =  \{1,\cdots,n-N\},
	$
    	$
    	k = \{1,\cdots,N \},
	$
and	
    \be
    \ba{rcl}
    \beta_{ik}
    & = &
    {\bf card} \left( \mre_i^T \! \left( I_s \, \circ \, K_r \right) \circ v_k^T \right)
    \\[0.15cm]
    \beta_i
    & = &
    {\bf card} \left( \mre_i^T \! \left( I_s \, \circ \, K_r \right) \right).
    \ea
    \label{eq.gr4}
    \ee
    \end{subequations}
The elementwise penalty~\eqref{eq.gr1} eliminates individual components of the feedback gain. In contrast, the group penalties~\eqref{eq.gr2} and~\eqref{eq.gr3} simultaneously eliminate feedback gains associated with a particular generator or feedback gains associated with all other generators, respectively. The {\em cardinality} function ${\bf card}(\cdot)$ in~\eqref{eq.gr4} counts the number of nonzero elements of a matrix, $\circ$ is elementwise matrix multiplication, $I_s \in \bbR^{m \times (n-N)}$ is the structural identity matrix (see Fig.~\ref{fig.IS} for the structure of $I_s$), $\mre_i \in \bbR^m$ is the $i$th unit vector, and $v_k \in \bbR^{n-N}$ is the structural identity vector. This vector is partitioned conformably with the partition of the vector $r$,
	\[
    v_k
    \; \DefinedAs \;
    \obth{\vartheta_1^T}{\cdots}{\vartheta_N^T}^T
    \]%
where $\vartheta_l = \one$ for $l=k$ and $\vartheta_l = 0$ for $ l \neq k$.

We note that the Euclidean norm ($\| \, \cdot \, \|_2$, not its square) is a widely used regularizer for enhancing group sparsity~\cite{yuan2006model}. The group weights $W_{i k}$'s and $W_i$'s are updated iteratively using the solution to~\eqref{eq.ADMM2} from the previous iteration~\cite{canwakboy08}. The scaling factors $\beta_{ik}$ and $\beta_i$ account for variations in the group sizes.

	\vspace*{-1ex}
\section{Case study: IEEE $39$ New England model}
\label{Section:Example}

\begin{figure}[!htb]
\centering
\begin{tabular}{c}
\includegraphics[width=0.4\textwidth]{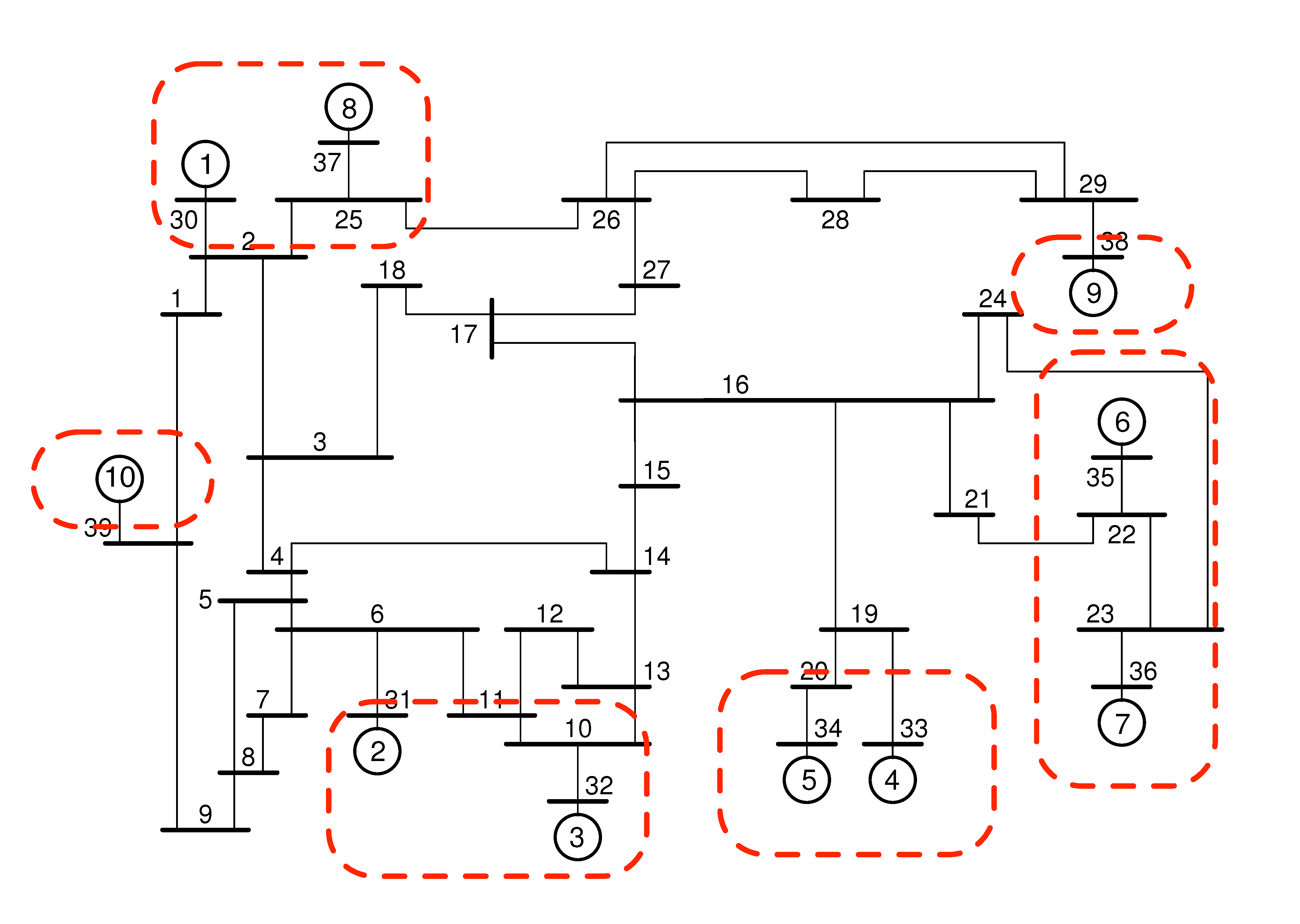}
\end{tabular}
\caption{The IEEE $39$ New England Power Grid and its coherent groups identified using slow coherency theory.}
\label{fig.NewEngland}
\end{figure}

The IEEE 39 New England Power Grid model consists of $39$ buses and $10$ detailed two-axis generator models; see Fig.~\ref{fig.NewEngland}. All loads are modeled as constant power loads. Generators $1$ to $9$ are equipped with PSSs, and generator $10$ is an equivalent aggregated model representing the transmission network of a neighboring area. This generator has an inertia which is an order of magnitude larger than the inertia of other generators.

The uncontrolled open-loop system is unstable, and PSSs are used for stabilization and to suppress local oscillations. For the subsequent analysis and the wide-area control design, we assume that {\em the PSS inputs are embedded in the open-loop matrix\/} $A \in \bbR^{75 \times 75}$ in~\eqref{eq.ss}. The transfer function of the local PSS controller on the $i$th generator is given by
	\be
    u_i (s)
    \; = \;
    k_i \cdot \displaystyle{\dfrac{T_{w,i} s}{1+T_{w,i} s}}
     \cdot \displaystyle{\dfrac{1+T_{n1,i} s}{1+T_{d1,i} s}}
     \cdot \displaystyle{\dfrac{1+T_{n2,i} s}{1+T_{d2,i} s}}
     \cdot \dot{\theta}_i (s)
     \non
    \ee
with controller gains $T_{w,i} = 5$, $T_{n1,i} = T_{n2,i} = 0.1$, $T_{d1,i} = T_{d2,i} = 0.01$, $k_i = 3$ for $i \in \{ 1, \cdots , 9 \}$. This set of PSS control gains stabilizes the unstable open-loop system, but it still features several poorly-damped modes. Our objective is to augment the local PSS control strategy with an optimal wide-area controller in order to simultaneously guard against inter-area oscillations and weakly dampened local oscillations.

Our computational experiments can be reproduced using the code available at:
	\begin{center}
{\small \mbox{www.umn.edu/$\sim$mihailo/software/lqrsp/matlab-files/lqrsp{\textunderscore}wac.zip}}
	\end{center}

	\vspace*{-3ex}
\subsection{Analysis of the open-loop system}
	\label{sec.ol-analysis}

\begin{figure}[t]
\centering
\begin{tabular}{ccc}
\subfloat[Mode $1$]
{\includegraphics[width=0.15\textwidth]{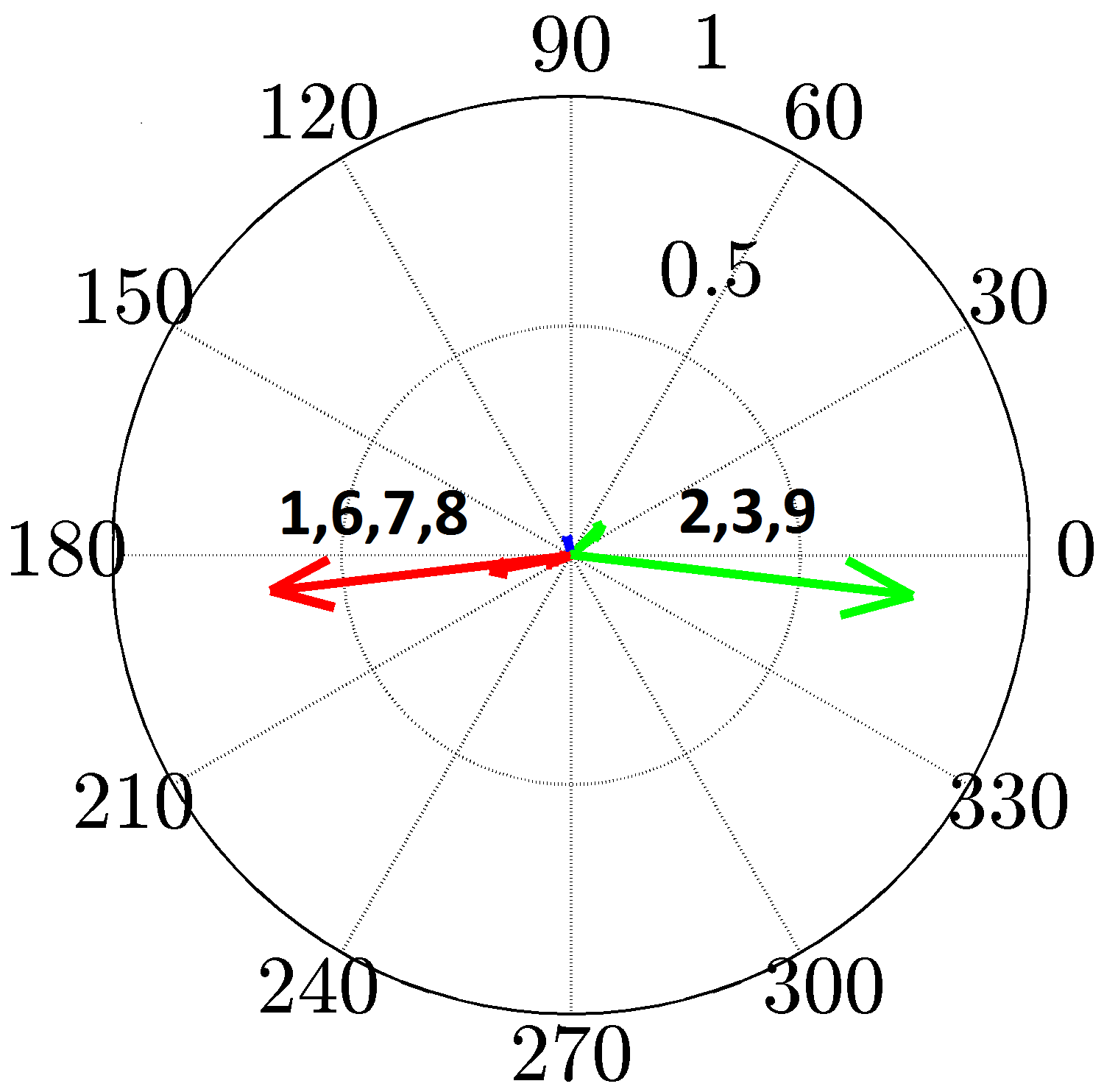}
\label{fig.mode1}}
\!\!\!\! & \!\!\!\!
\subfloat[Mode $2$]
{\includegraphics[width=0.15\textwidth]{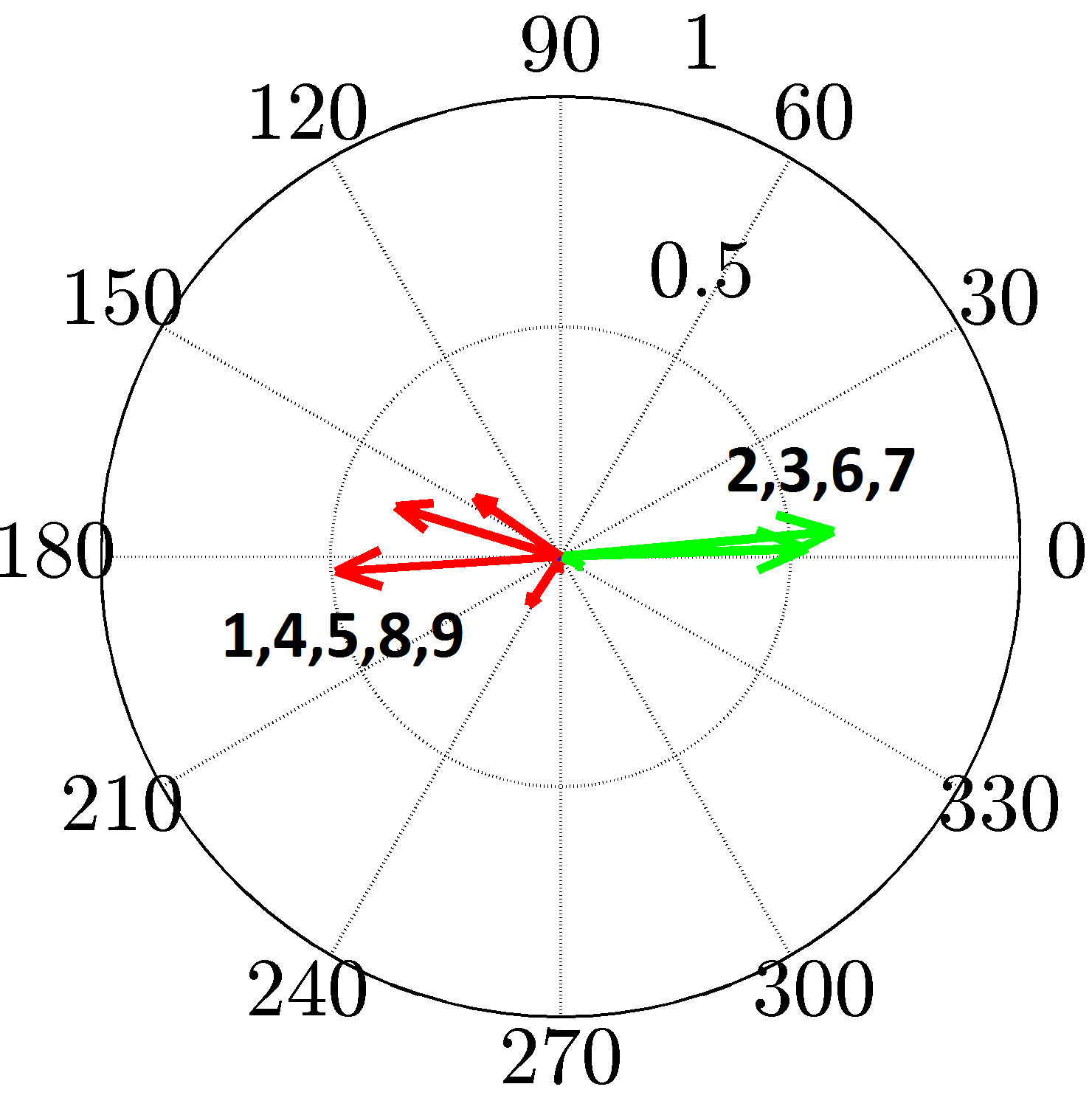}
\label{fig.mode2}}
\!\!\!\! & \!\!\!\!
\subfloat[Mode $3$]
{\includegraphics[width=0.15\textwidth]{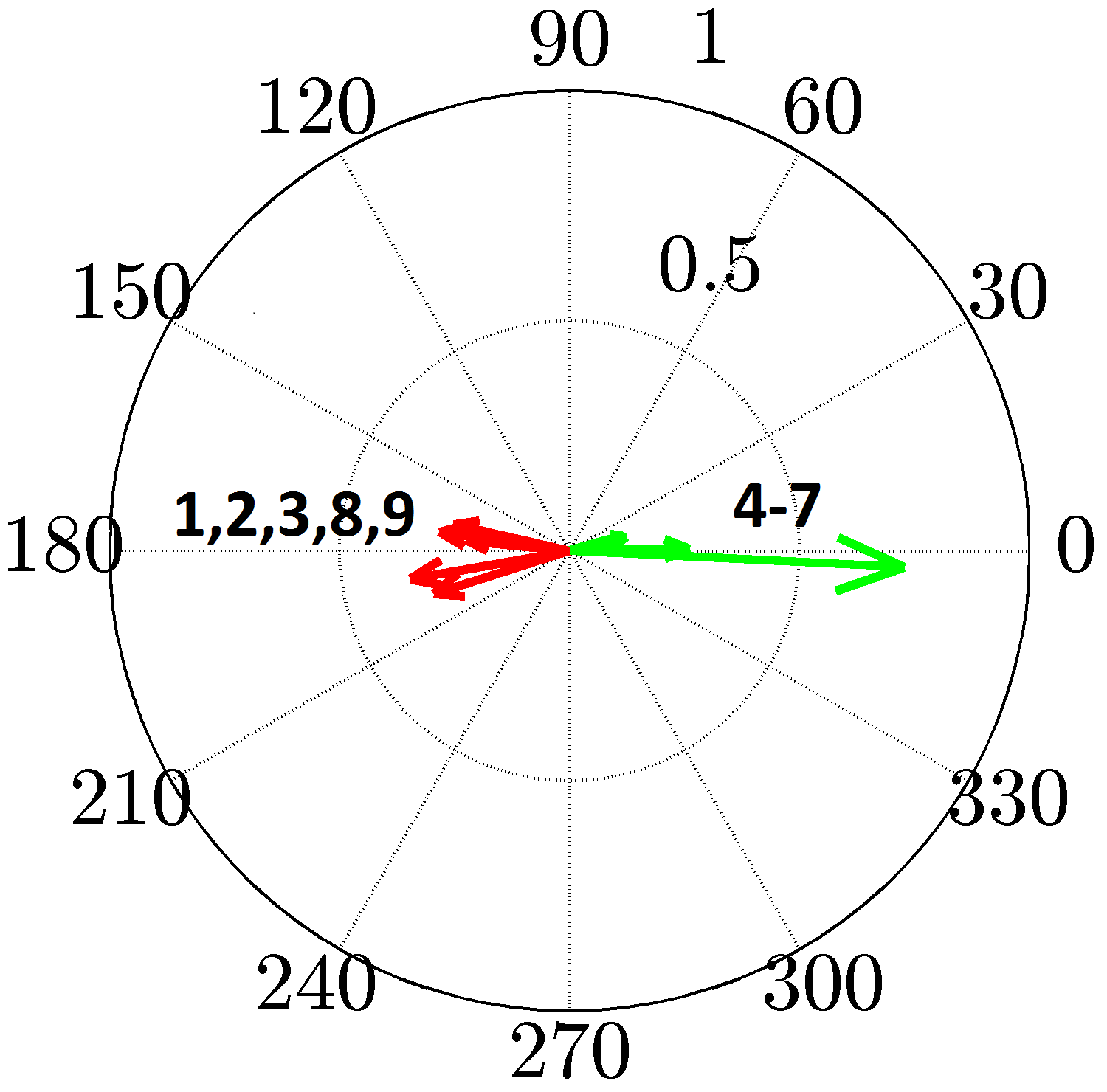}
\label{fig. mode3}}
\\
\subfloat[Mode $4$]
{\includegraphics[width=0.15\textwidth]{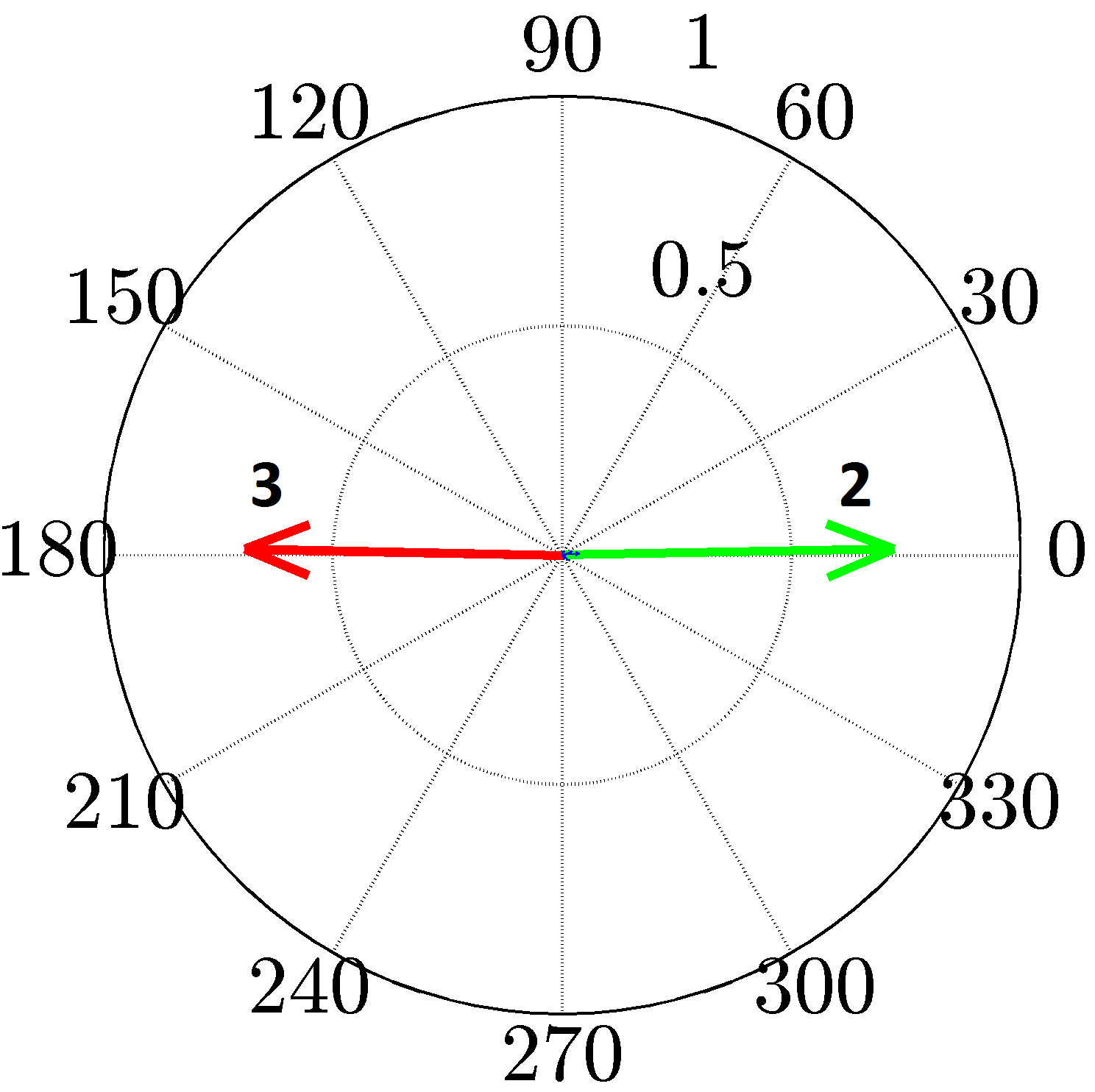}
\label{fig.mode4}}
\!\!\!\! & \!\!\!\!
\subfloat[Mode $5$]
{\includegraphics[width=0.15\textwidth]{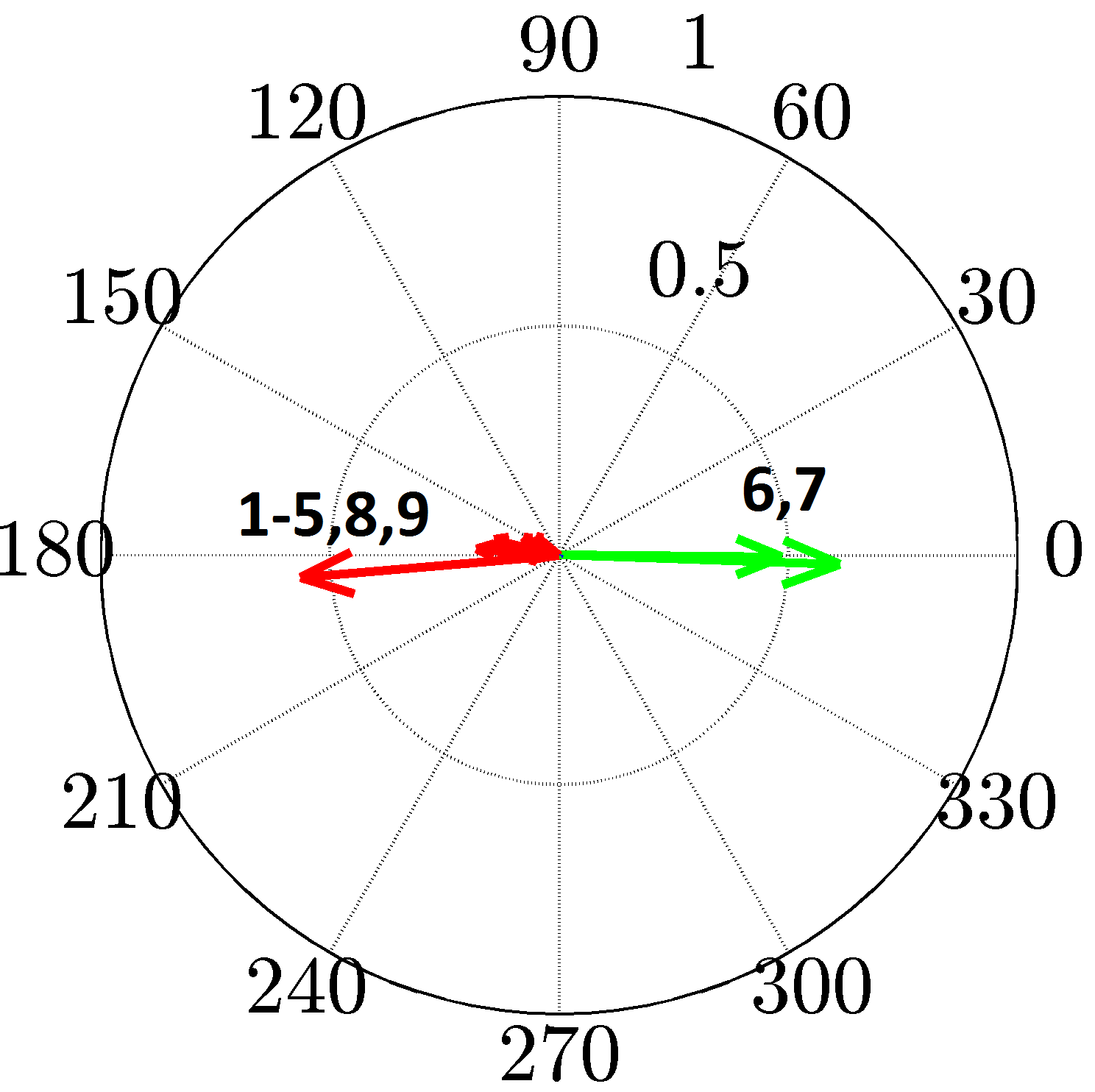}
\label{fig. mode5}}
\!\!\!\! & \!\!\!\!
\subfloat[Mode $6$]
{\includegraphics[width=0.15\textwidth]{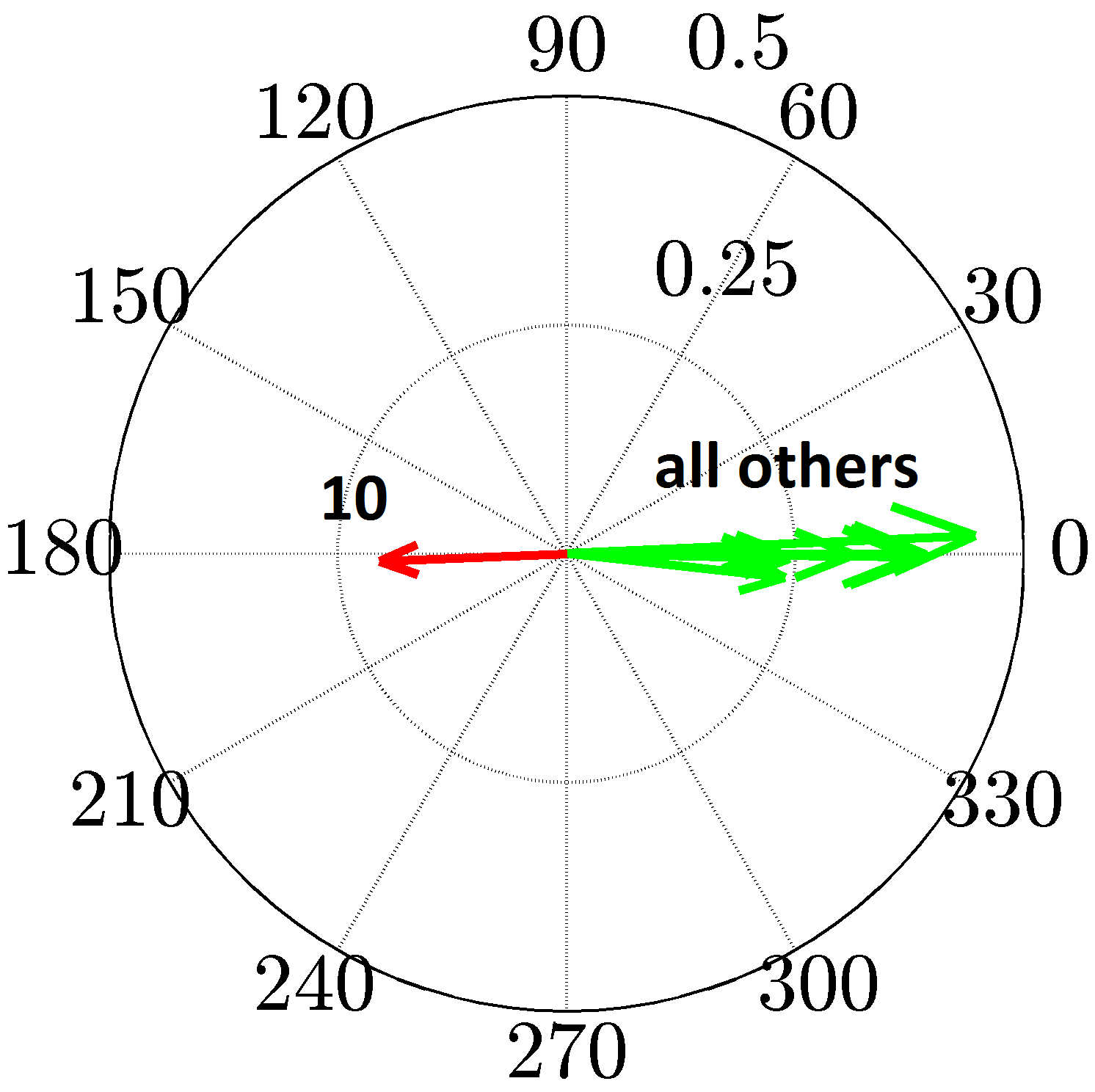}
\label{fig. mode6}}
\end{tabular}
\caption{Polar plots of the angle components of the six poorly-damped modes for the open-loop system.}
\label{fig.modes}
\vspace*{0.25cm}
\end{figure}

Despite the action of the local PSS controllers, modal and participation factor analyses reveal the presence of six poorly-damped modes in the New England power grid model; see Table~\ref{Table: poorly-damped modes of New England Grid} and Fig.~\ref{fig.modes}. Mode $4$ is a local mode because it only involves oscillations between generators $2$ and $3$, which belong to the same coherent group. All other modes are inter-area modes where groups of generators oscillate against each other. Since these inter-area modes are  poorly damped with damping ratios as low as $1.20 \%$ and $2.61 \%$, the local PSS controllers need to be complemented by  supplementary wide-area controllers to improve the damping of the inter-area oscillations.

We depart from the modal perspective and examine the power spectral density and variance amplification of the open-loop system. This type of analysis allows us to identify (i) the temporal frequencies for which large amplification occurs; and (ii) the spatial structure of strongly amplified responses.

\begin{table}[t]
\small
\centering
\caption{Poorly-damped modes of New England model}
\label{Table: poorly-damped modes of New England Grid}
\begin{tabular}{| c | l | l | l | l |}
\hline
\!\!\!mode\,\!\! \!\!\!\!&\!\!\!\! eigenvalue \!\!\!\!&\!\!\!\! damping \!\!\!\!\!&\!\!\!\! freq. \!\!\!\!\!\!&\!\!\!\! coherent \\
no. \!\!\!\!&\!\!\!\! pair  \!\!\!\!&\!\!\!\! ratio \!\!\!\!\!&\!\!\!\! $[\textup{Hz}]$ \!\!\!\!\!&\!\!\!\! groups \\
\hline\hline
\!\!\!1 \!\!\!\!&\!\!\!\! $-0.0882 \pm \textup{j}\, 7.3695$ \!\!\!\!&\!\!\!\! 0.0120 \!\!\!\!\!&\!\!\!\! 1.1618 \!\!\!\!\!&\!\!\!\! 1,6,7,8 vs. 2,3,9 \!\!\!\!\!\\
\!\!\!2 \!\!\!\!&\!\!\!\! $-0.1788 \pm \textup{j}\, 6.8611$ \!\!\!\!&\!\!\!\! 0.0261 \!\!\!\!\!&\!\!\!\! 1.0918 \!\!\!\!\!&\!\!\!\! 2,3,6,7 vs.1,4,5,8,9 \!\!\!\!\!\\
\!\!\!3 \!\!\!\!&\!\!\!\! $-0.2404 \pm \textup{j}\, 6.5202$ \!\!\!\!&\!\!\!\! 0.0368 \!\!\!\!\!&\!\!\!\! 1.0377 \!\!\!\!\!&\!\!\!\! 1,2,3,8,9 vs. 4-7 \!\!\!\!\!\\
\!\!\!4 \!\!\!\!&\!\!\!\! $-0.4933 \pm \textup{j}\, 7.7294$ \!\!\!\!&\!\!\!\! 0.0637 \!\!\!\!\!&\!\!\!\! 1.2335 \!\!\!\!\!&\!\!\!\! 2 vs. 3 \!\!\!\!\!\\
\!\!\!5 \!\!\!\!&\!\!\!\! $-0.4773 \pm \textup{j}\, 6.9858$ \!\!\!\!&\!\!\!\! 0.0682 \!\!\!\!\!&\!\!\!\! 1.1141 \!\!\!\!\!&\!\!\!\! 6,7 vs. 1-5,8,9 \!\!\!\!\!\\
\!\!\!6 \!\!\!\!&\!\!\!\! $-0.3189 \pm \textup{j}\, 4.0906$ \!\!\!\!&\!\!\!\! 0.0777 \!\!\!\!\!&\!\!\!\! 0.6525 \!\!\!\!\!&\!\!\!\! 10 vs. all others \!\!\!\!\!\\
\hline
\end{tabular}
\end{table}

\begin{figure}[!htb]
\centering
\begin{tabular}{cc}
\subfloat[]
{\includegraphics[width=0.22\textwidth]{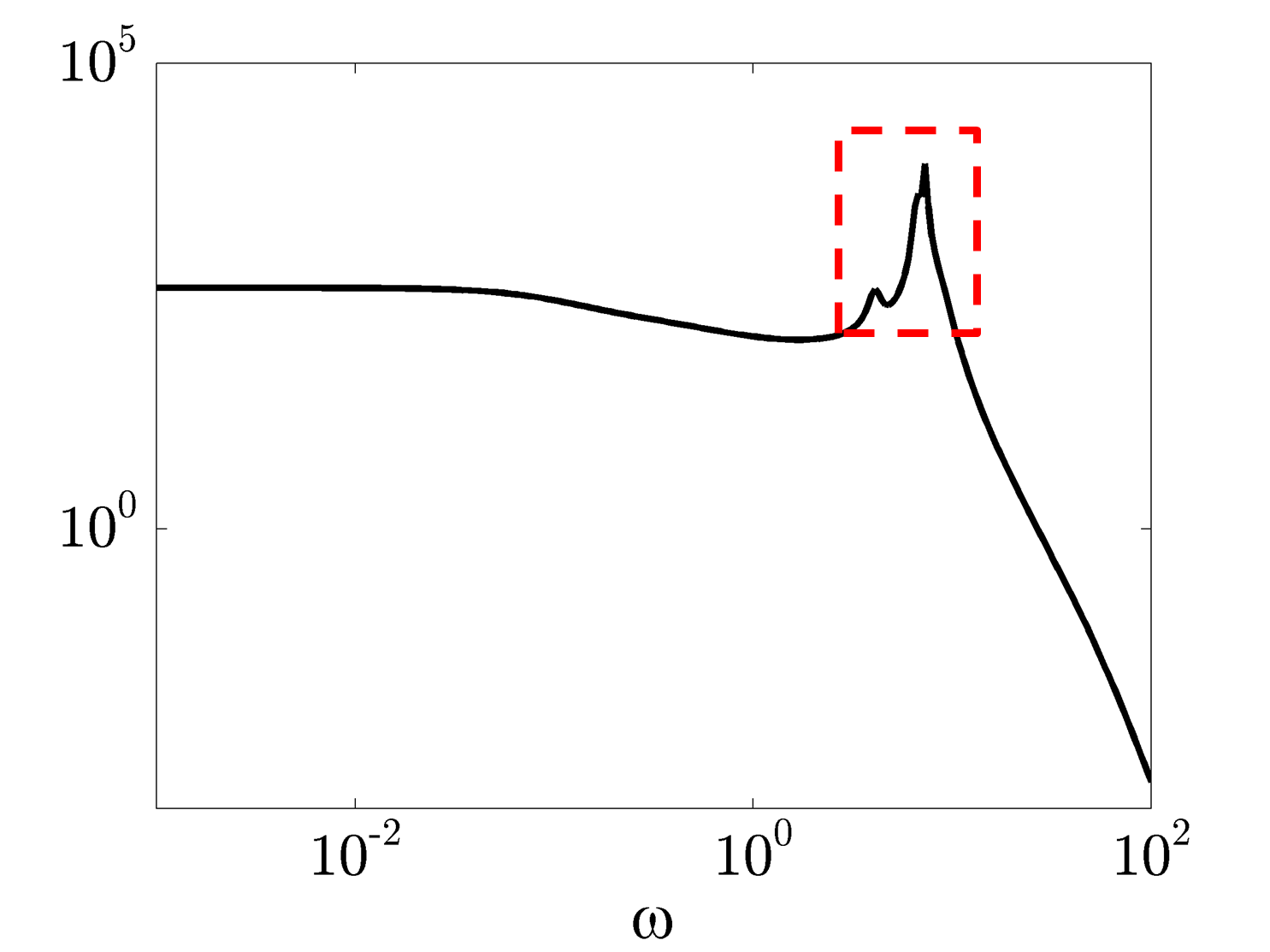}
 \label{fig.psd_pss5}}
&
\subfloat[]
{\includegraphics[width=0.22\textwidth]{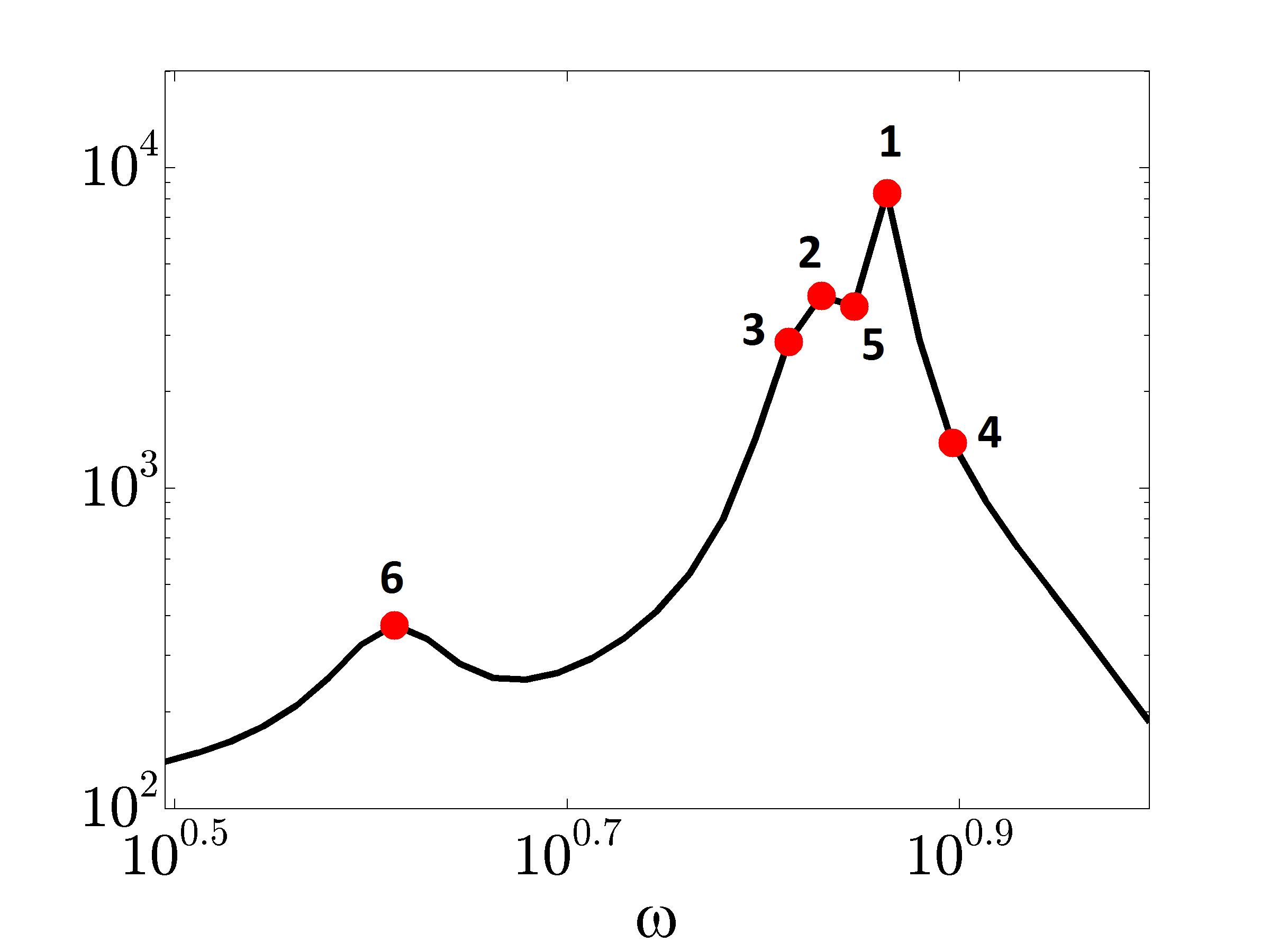}
 \label{fig.psd_pss5_zoom}}
\end{tabular}
\vspace*{0.1cm}
\caption{(a) Power spectral density of the open-loop system; (b) zoomed version of the red square shown in (a). Red dots denote poorly-damped modes from Table~\ref{Table: poorly-damped modes of New England Grid}.}
\label{fig.psd_open}
\end{figure}

Figure~\ref{fig.psd_open} illustrates the power spectral density of the open-loop system. The largest peak occurs at $\omega_1 = 7.2925$ rad/s ($f_1 = {\omega_1}/{2 \pi} = 1.1606$ Hz) and it corresponds to mode $1$ in Table~\ref{Table: poorly-damped modes of New England Grid} and Fig.~\ref{fig.modes}. Another resonant peak at $\omega_2 = 4.0930$ rad/s ($f_2 = 0.6514$ Hz) corresponds to mode $6$ in Table~\ref{Table: poorly-damped modes of New England Grid} and Fig.~\ref{fig.modes}. The red dots in Fig.~\ref{fig.psd_pss5_zoom} indicate all six poorly-damped modes.

\begin{figure}[]
\centering
\begin{tabular}{c}
\includegraphics[width=0.4\textwidth,height=0.2\textwidth]{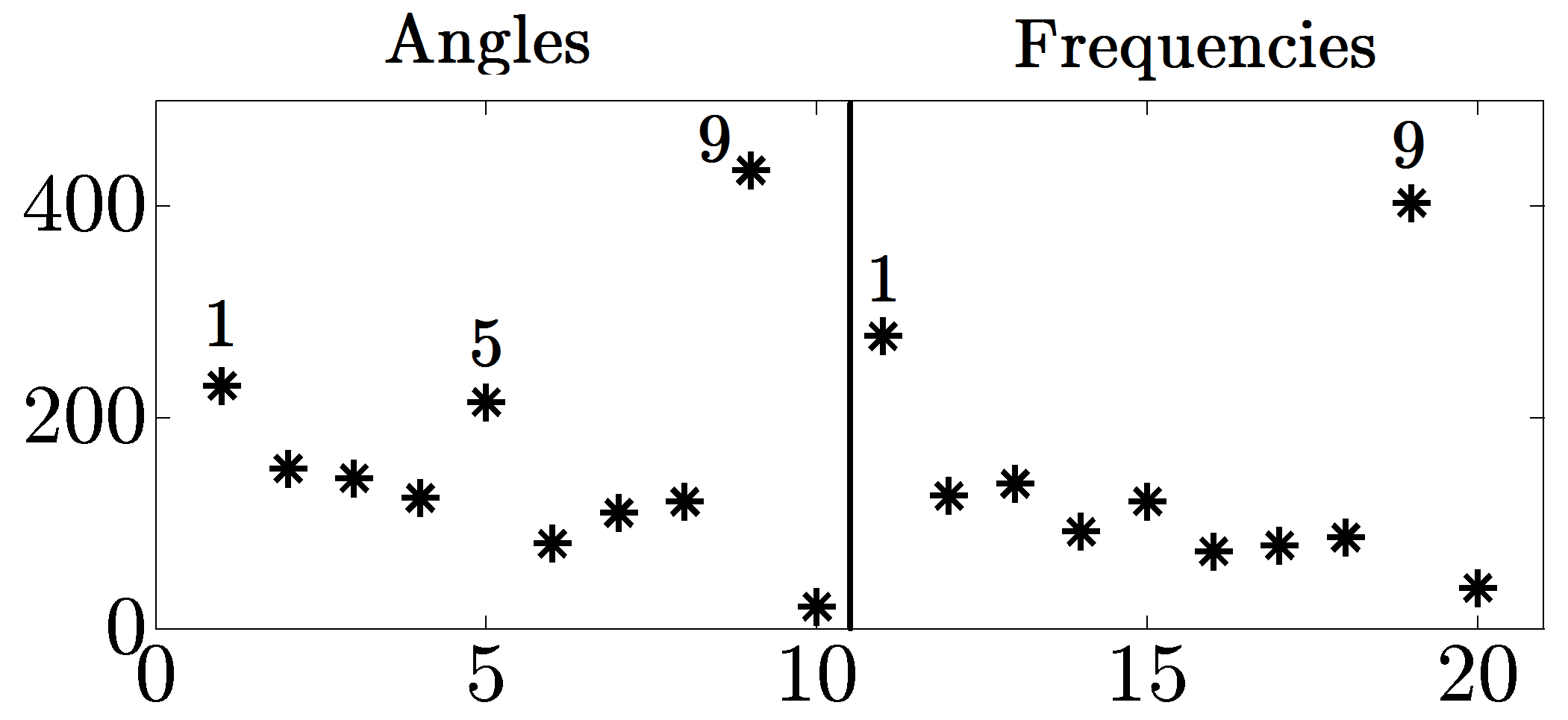}
\end{tabular}
\caption{Diagonal elements of the open-loop covariance matrix $Z_1$ determine contribution of each generator to the variance amplification.}
\label{fig.va_open}
\end{figure}

The contribution of each generator to the steady-state variance is shown in Fig.~\ref{fig.va_open}. The diagonal elements of the output covariance matrix $Z_1$ contain information about mean-square deviation from angle average and variance amplification of frequencies of the individual generators. From Fig.~\ref{fig.va_open}, we see that the largest contribution to the variance amplification arises from the misalignment of angles of generators $1$, $5$, and $9$, and misalignment of frequencies of generators $1$ and $9$.

\begin{figure}[]
\centering
\begin{tabular}{cc}
\subfloat[Eigenvalues of $Z_1$]
{\includegraphics[width=0.22\textwidth]{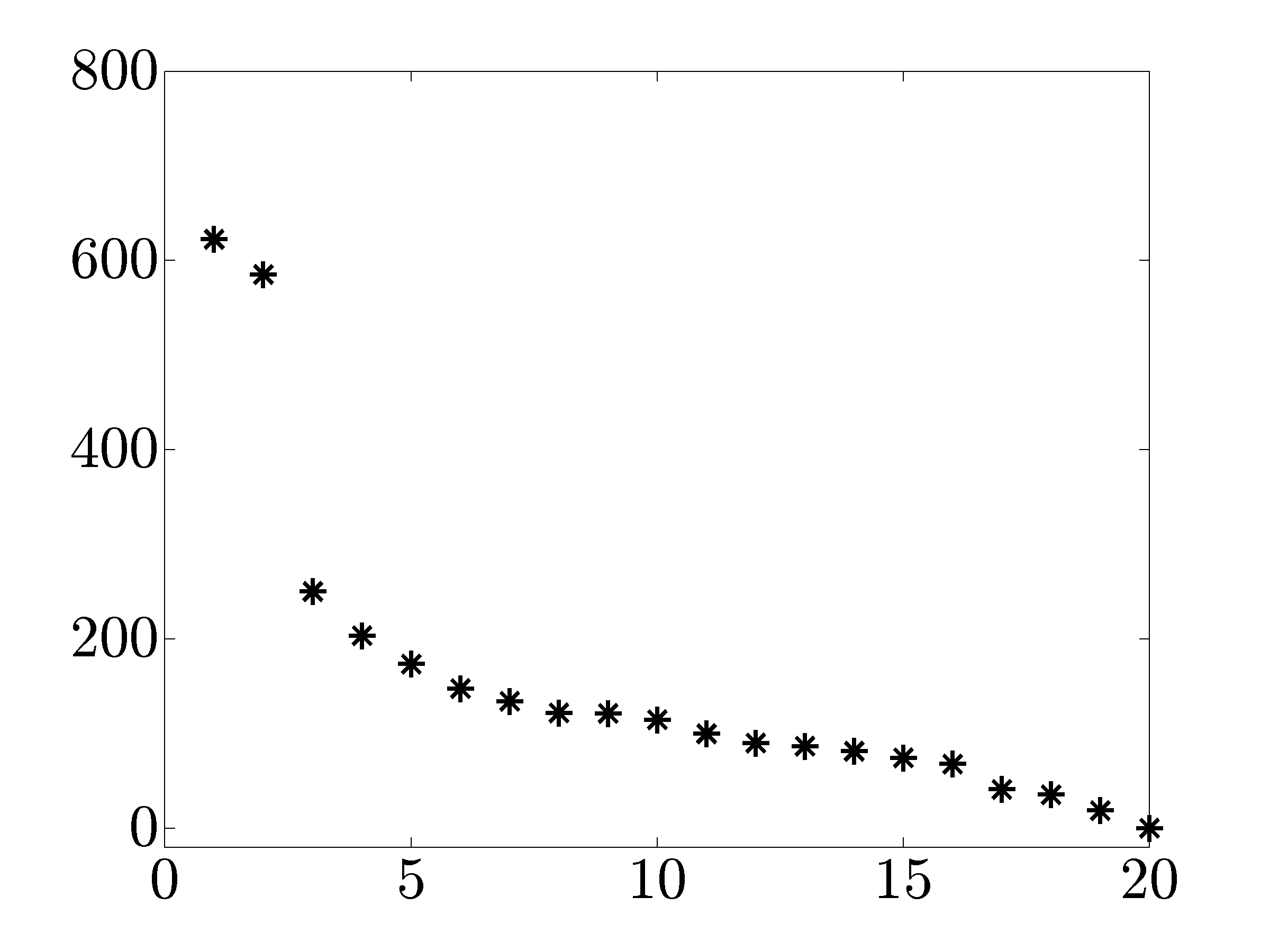}
\label{fig.eig_open}}
&
\subfloat[$\lambda_1 (Z_1)$]
{\includegraphics[width=0.22\textwidth]{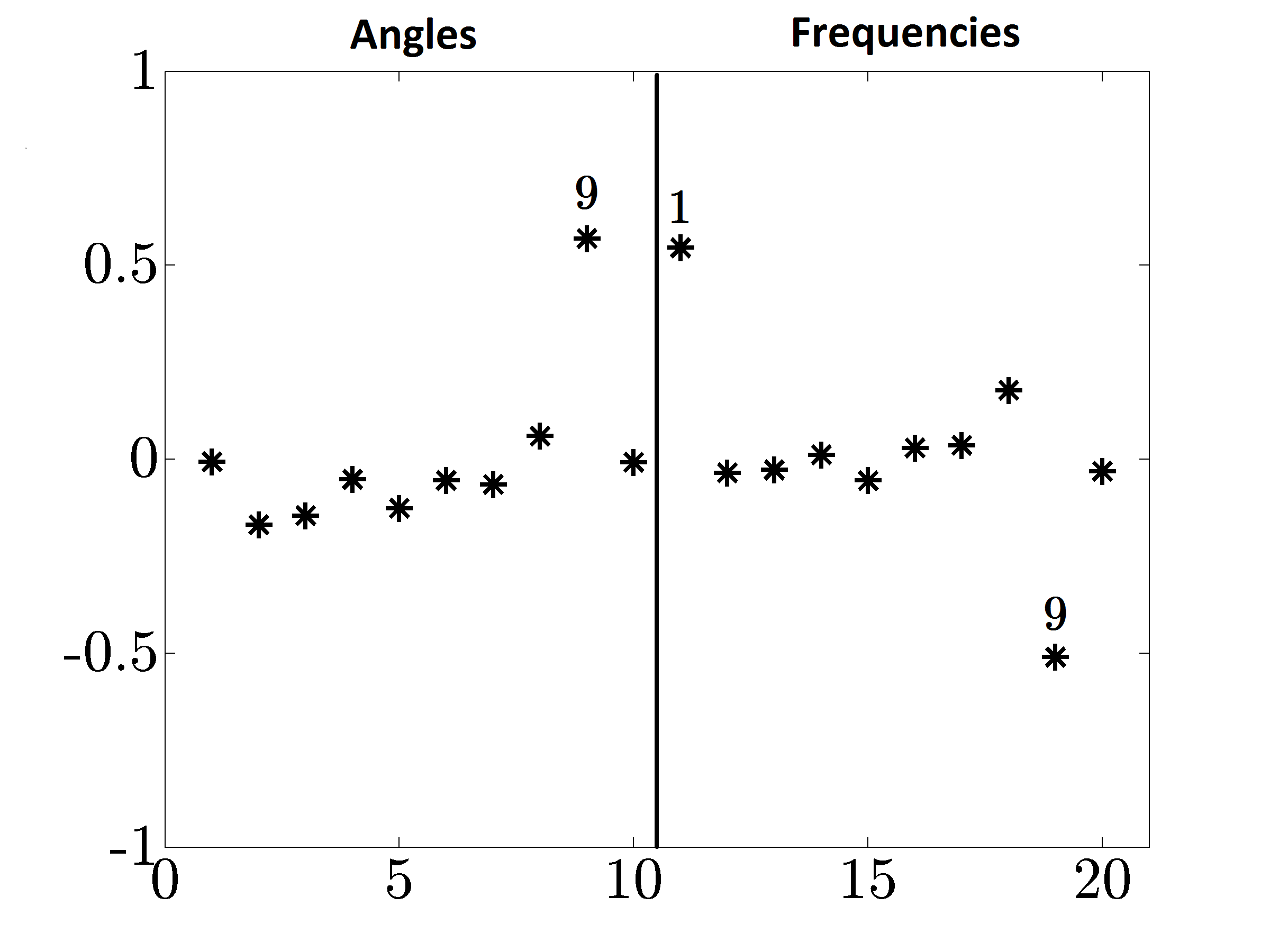}} \\
\subfloat[$\lambda_2 (Z_1)$]
{\includegraphics[width=0.22\textwidth]{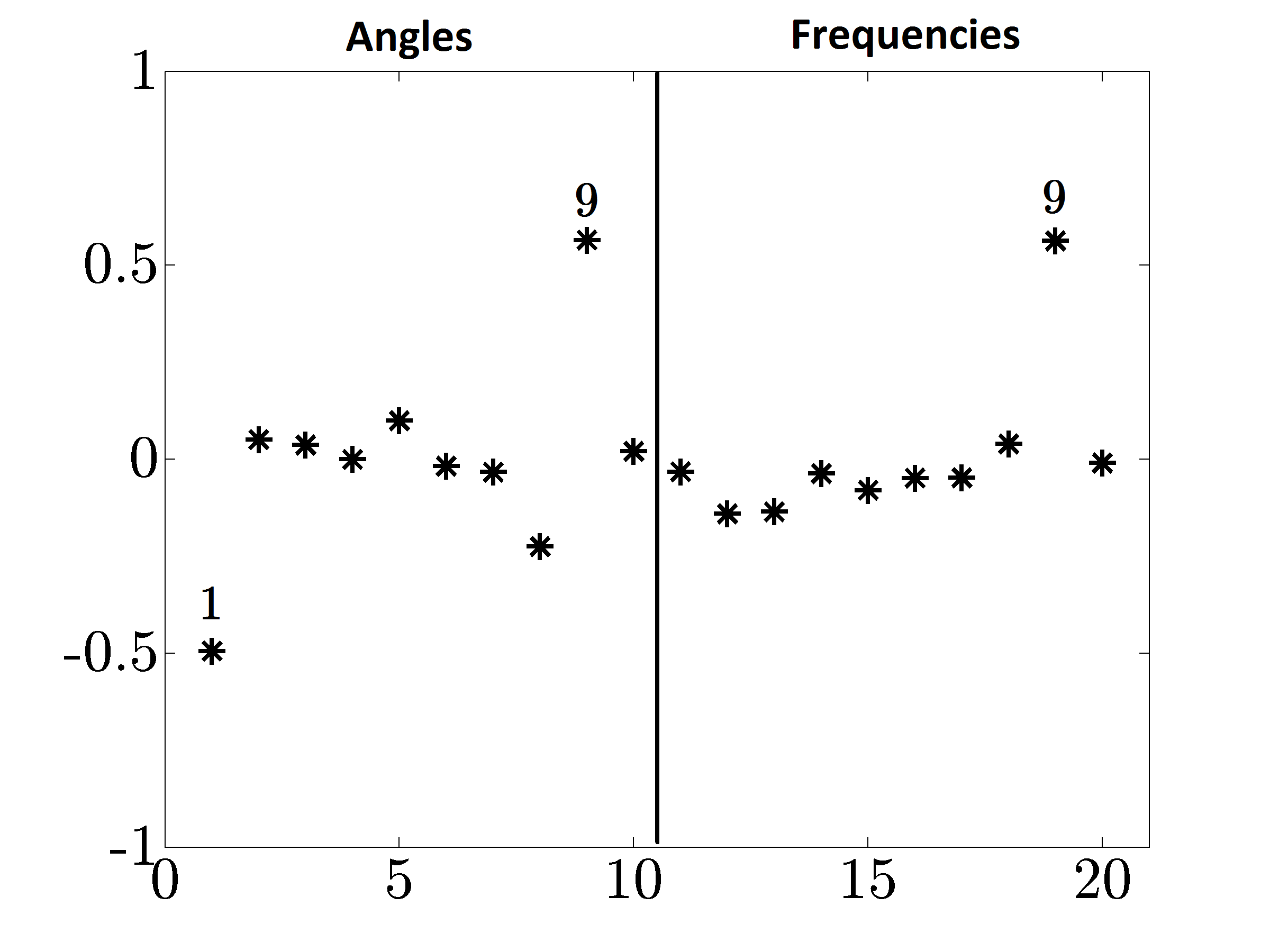}}
&
\subfloat[$\lambda_3 (Z_1)$]
{\includegraphics[width=0.22\textwidth]{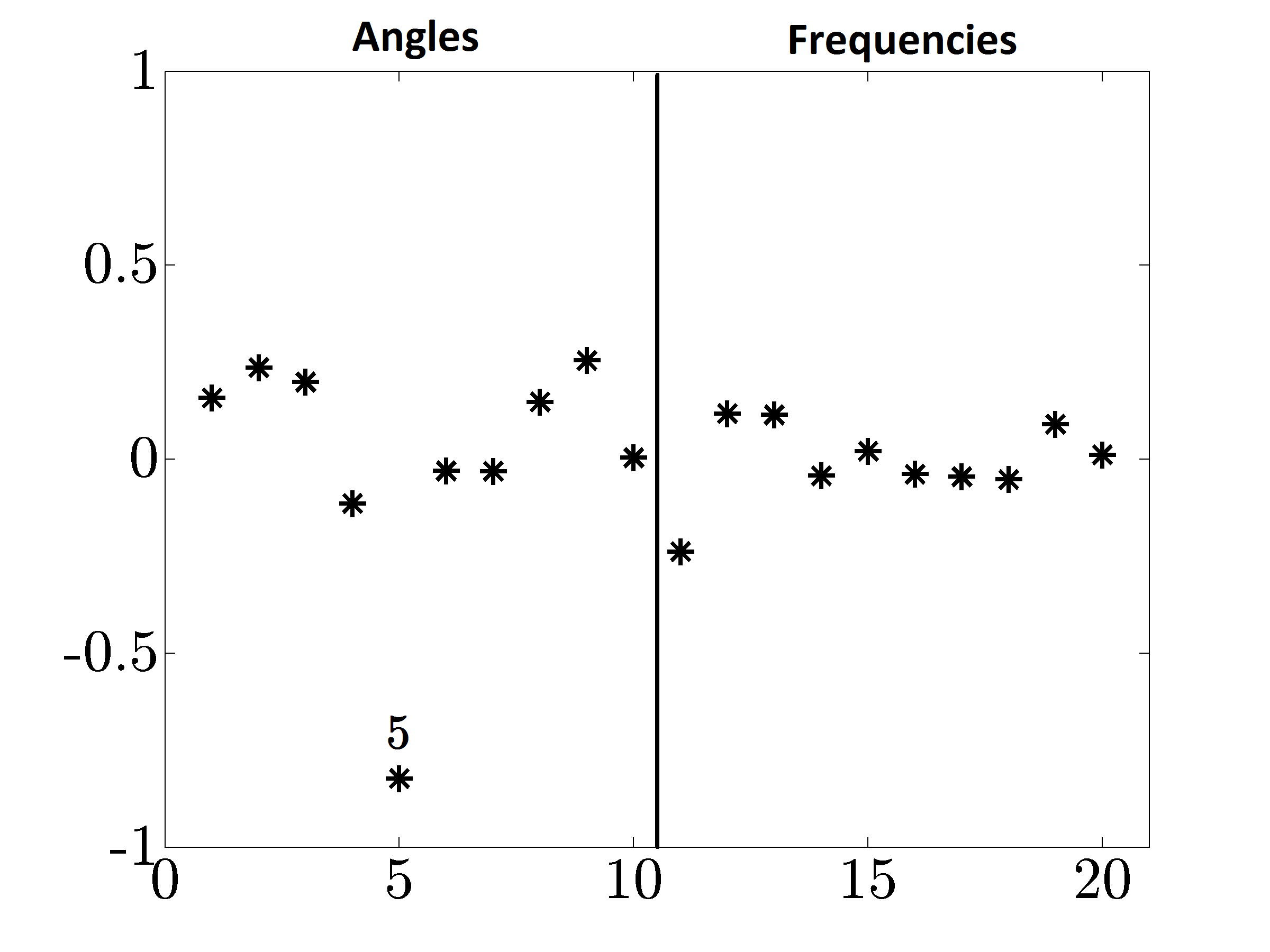}}
\end{tabular}
\caption{(a) Eigenvalues; and (b)-(d) eigenvectors corresponding to the three largest eigenvalues $\lambda_i$ of the open-loop output covariance matrix $Z_1$.}
\label{fig.maxeig}
\end{figure}

Similar observations can be made from Fig.~\ref{fig.maxeig}. In Fig.~\ref{fig.eig_open}, we observe two dominant eigenvalues of the output covariance matrix $Z_1$. We also show the spatial structure of the three principal eigenvectors (modes) of $Z_1$, which contain $47.5 \%$ of the total variance. Although the angle and frequency fluctuations in experiments and nonlinear simulations are expected to be more complex than the structures presented in Fig.~\ref{fig.maxeig}, the spatial profiles identified here are likely to play significant role in amplification of disturbances in power systems.

	\vspace*{-2ex}
\subsection{Sparsity-promoting optimal wide-area control}
\label{Section:sparsity-promoting}

We next illustrate that the addition of certain long-range communication links and careful retuning of the local excitation controllers are effective means for improving the system performance and increasing its resilience to inter-area oscillations.

\subsubsection{Elementwise sparsity}
\label{Section:element}

\begin{figure}[!htb]
\centering
\subfloat[$\gamma = 0.0818$, $\card \, (K) = 43$]
	{
\includegraphics[width=0.45\textwidth]{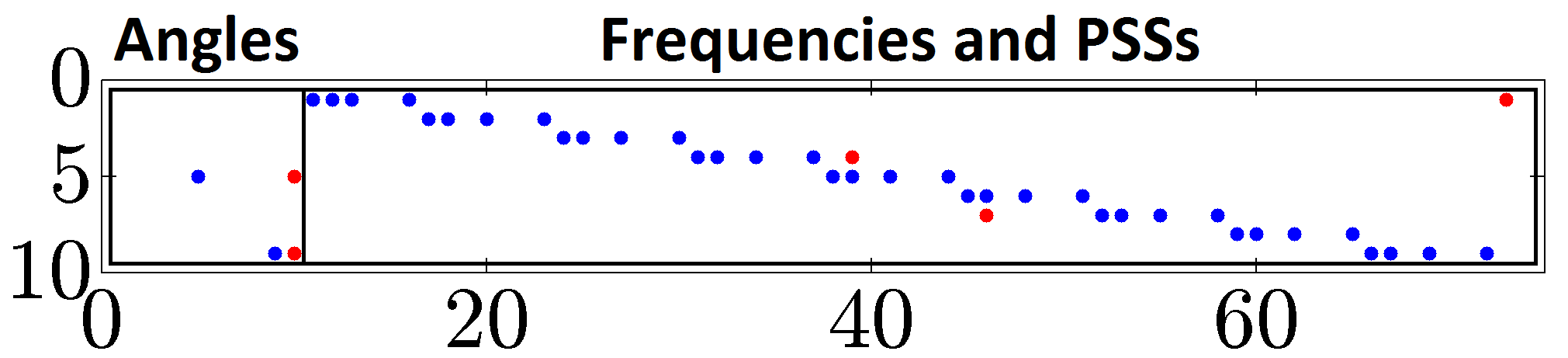}
\label{fig.Ga_M}
	}
	\\[-0.02cm]
\subfloat[$\gamma = 0.1548$, $\card \, (K) = 38$]
	{
\includegraphics[width=0.45\textwidth]{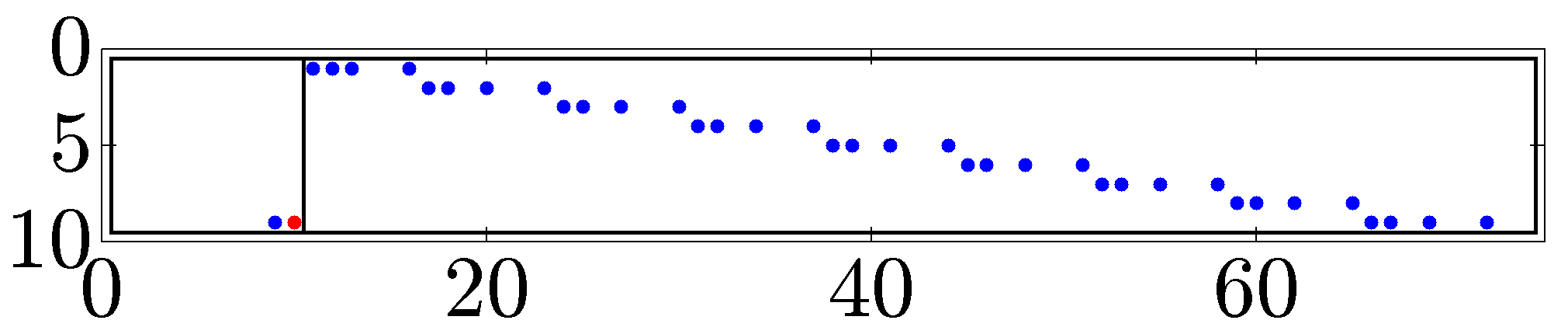}
\label{fig.Gb_M}
	}
	\\[-0.02cm]
\subfloat[$\gamma = 0.2500$, $\card \, (K) = 35$]
	{
\includegraphics[width=0.45\textwidth]{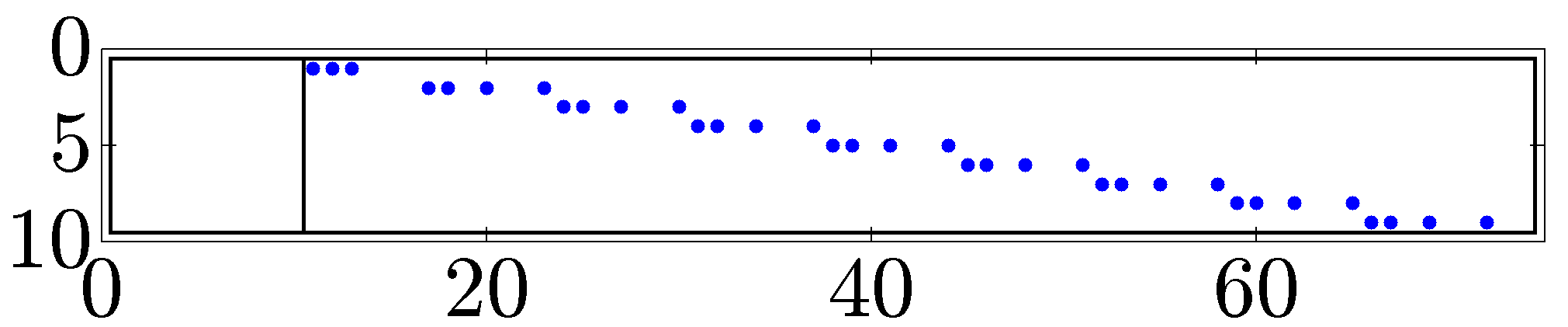}
\label{fig.Gc_M}
}
\caption{Sparsity patterns of $K$ resulting from~\eqref{eq.ADMM1}.}
\label{fig.spy_coc}
\end{figure}

We first consider an optimal sparse controller whose structure is identified using the solution to~\eqref{eq.ADMM1}. Sparsity patterns of the feedback matrix $K \in \bbR^{9 \times 75}$ for different values of $\gamma$ are illustrated in Fig.~\ref{fig.spy_coc}. The blue dots denote information coming from the generators on which the particular controller acts, and the red dots identify information that needs to be communicated from other generators. For $\gamma = 0.0818$, the identified wide-area control architecture imposes the following requirements: (i) the controller of generator $9$, which contributes most to the variance amplification of both angles and frequencies, requires angle and field voltage measurements of the aggregate generator $10$; (ii) the controller of generator $5$ requires the difference between its angle and the angle of the equivalenced model $10$; and (iii) the controllers of generators $1$, $4$, and $7$ utilize the field voltage information of generators $10$, $5$, and $6$, respectively.

\begin{figure}[!htb]
\centering
\begin{tabular}{c}
\includegraphics[width=0.45\textwidth]{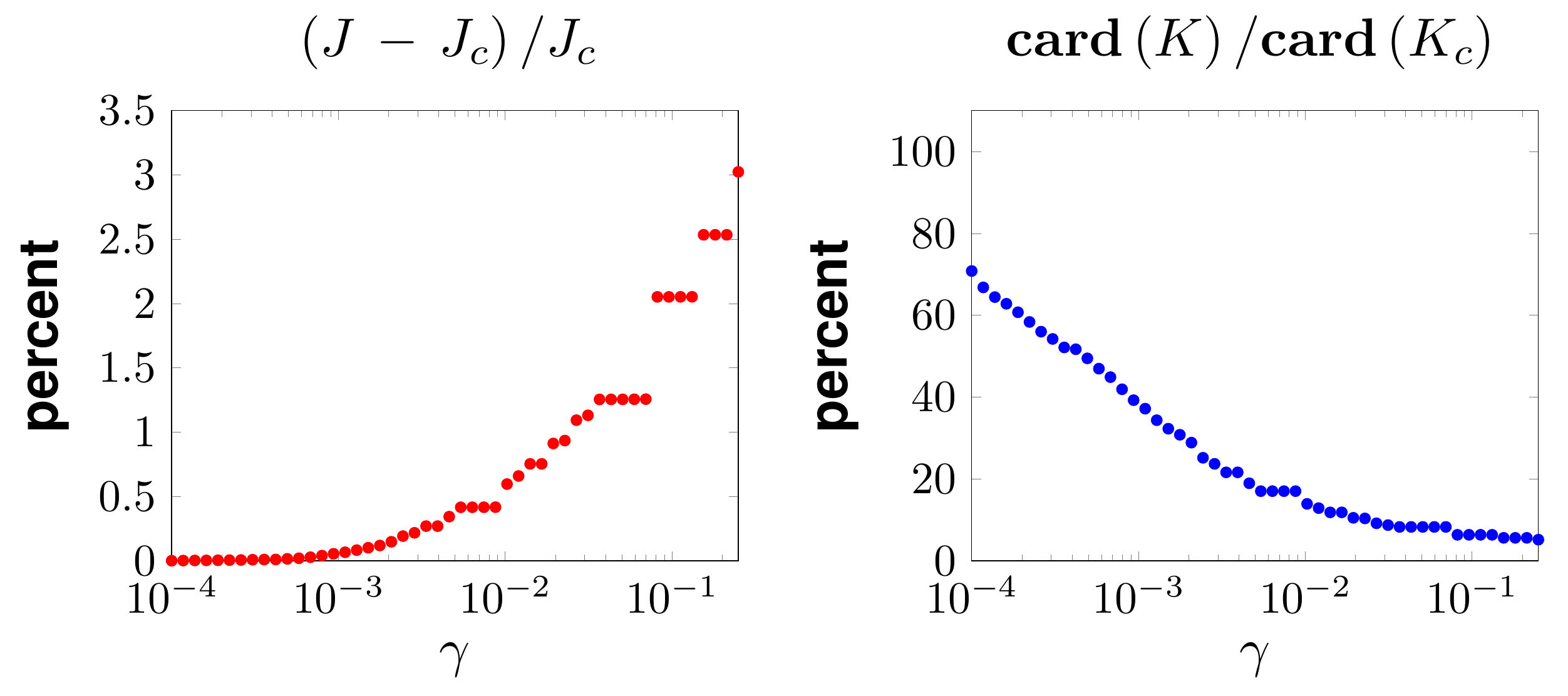}
\end{tabular}
\caption{Performance vs sparsity comparison of sparse $K$ and the optimal centralized controller $K_c$ for $50$ logarithmically-spaced points $ \gamma \in [ \, 10^{-4} \,,\, 0.25 \, ]$.}
\label{fig.sp_ne_coc}
\end{figure}

When $\gamma$ is increased to $0.1548$, only one long-range link remains. This link is identified by the red dot in Fig.~\ref{fig.Gb_M}, indicating that the controller of generator $9$ requires access to the angle mismatch relative to generator $10$. By further increasing $\gamma$ to $0.25$, we obtain a fully-decentralized controller. Compared to the optimal centralized controller, our fully-decentralized controller degrades the closed-loop performance by about $3.02\%$; see Fig.~\ref{fig.sp_ne_coc}. This fully-decentralized controller can be embedded into the local generator excitation system by directly feeding the local measurements to the automatic voltage regulator, thereby effectively retuning the PSS controller.

In earlier work~\cite{dorjovchebulACC13,dorjovchebulTPS14}, a small regularization term was added to the diagonal elements of the matrix $Q_{\theta}$ in order to provide detectability of the average mode. This has resulted in a controller that requires access to the absolute angle measurements to stabilize the average rotor angle. Our results indicate that long-range links identified in~\cite{dorjovchebulACC13,dorjovchebulTPS14} do not have significant influence on the system performance.

\subsubsection{Block sparsity}

\begin{figure}[!htb]
\centering
\subfloat[$\gamma = 0.0697$, $\card \, (K) = 66$]
	{
\includegraphics[width=0.45\textwidth]{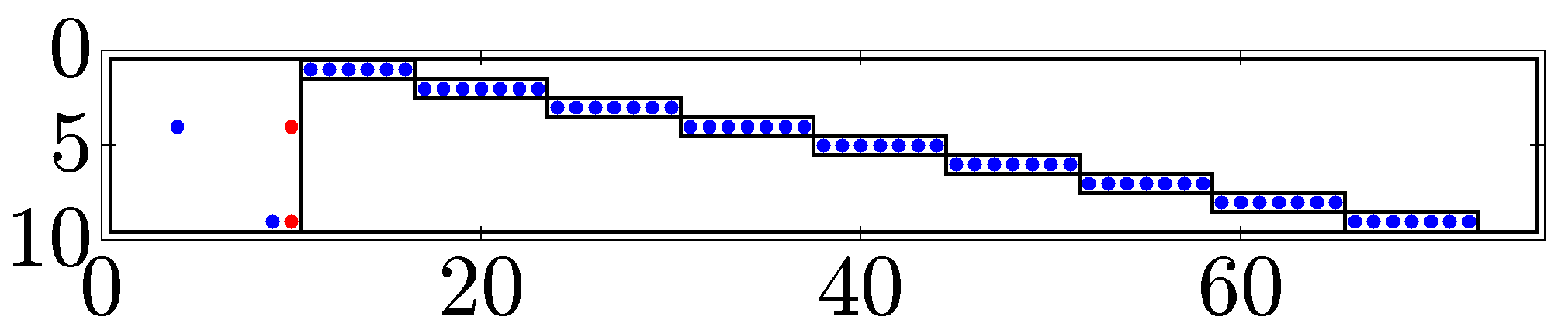}
\label{fig.G41}
	}
	\\[-0.02cm]
\subfloat[$\gamma = 0.0818$, $\card \, (K) = 64$]
	{
\includegraphics[width=0.45\textwidth]{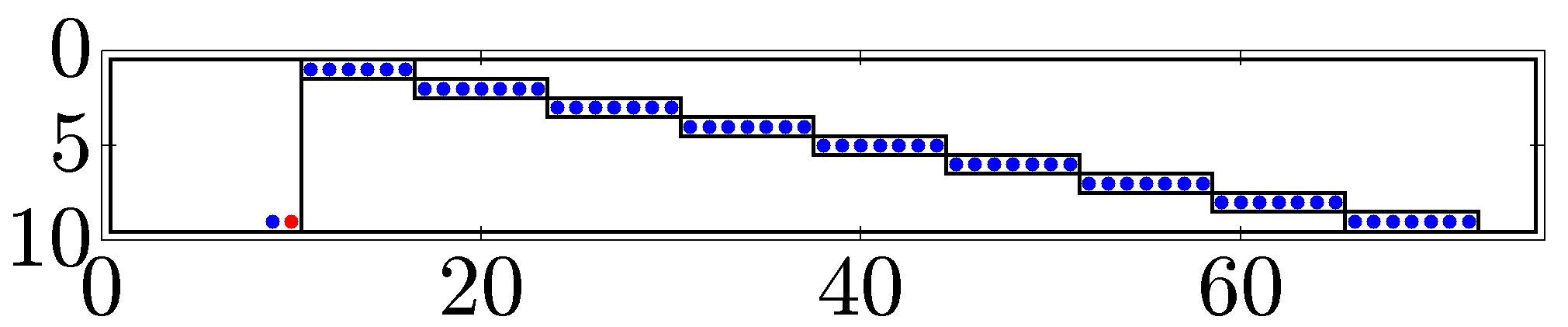}
\label{fig.G43}
	}
	\\[-0.02cm]
\subfloat[$\gamma = 0.2500$,  $\card \, (K) = 62$]
	{
\includegraphics[width=0.45\textwidth]{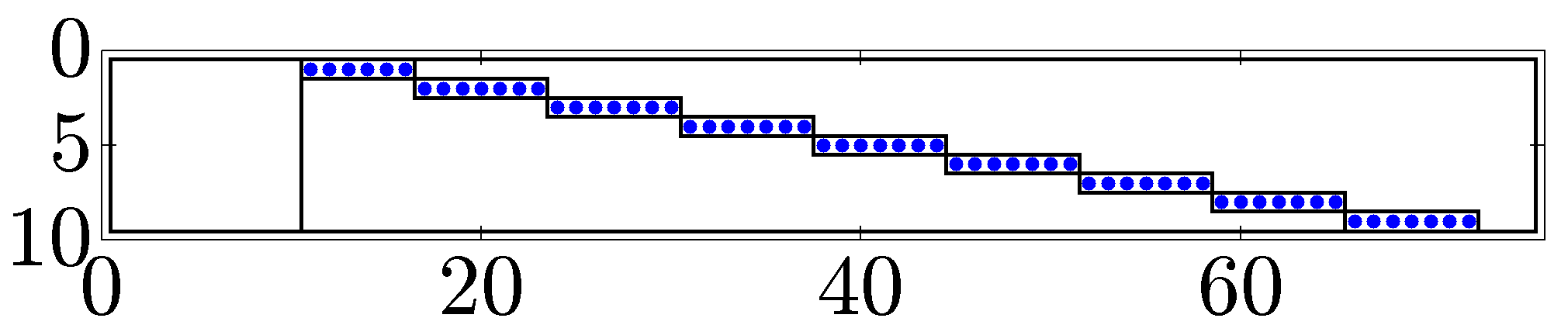}
\label{fig.G50}
}
\caption{Sparsity patterns of $K$ resulting from~\eqref{eq.ADMM2}.}
\label{fig.spy_blk}
\end{figure}

Three identified sparsity patterns of the feedback matrix resulting from the solution to~\eqref{eq.ADMM2}, with $g_\theta$ and $g_r$ given by~\eqref{eq.gtheta} and~\eqref{eq.gr3}, are shown in Fig.~\ref{fig.spy_blk}. In all three cases, structures of the angle feedback gains agree with the elementwise sparse controllers; cf.\ Fig.~\ref{fig.spy_coc}. On the other hand, the group penalty~\eqref{eq.gr3} yields block-diagonal feedback gains that act on the remaining states of generators $1$-$9$. Since no information exchange with aggregate generator $10$ is required, this part of the controller can be implemented in a fully-decentralized fashion in all three cases.

\begin{figure}[!htb]
\centering
\begin{tabular}{c}
\includegraphics[width=0.45\textwidth]{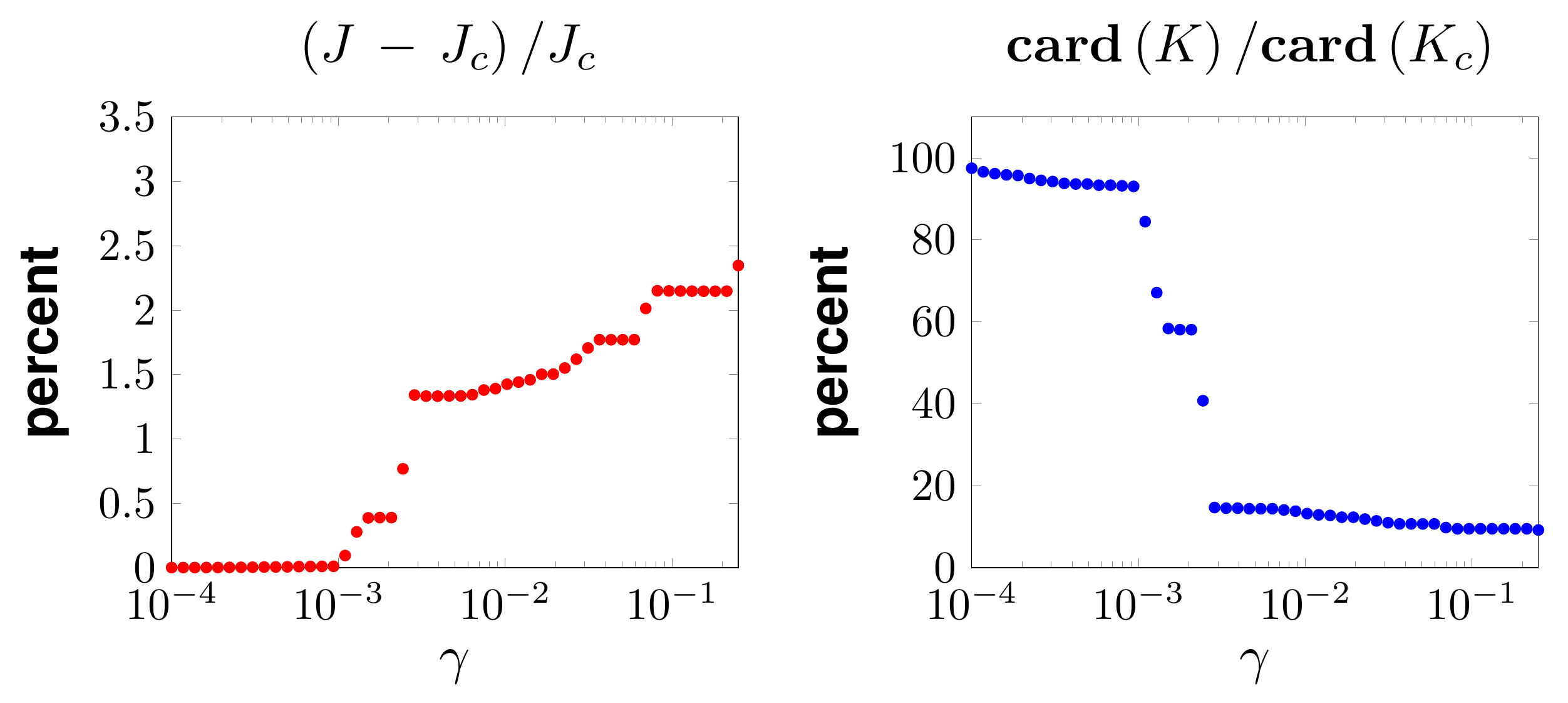}
\end{tabular}
\caption{Performance vs sparsity comparison of block-sparse $K$ and the optimal centralized controller $K_c$ for $50$ logarithmically-spaced points $ \gamma = \gamma_\theta = \gamma_r \in [ \, 10^{-4} \,,\, 0.25 \, ]$.}
\label{fig.perf_blk}
\end{figure}

Compared to the optimal centralized controller, a fully-decentralized controller with structure shown in Fig.~\ref{fig.G50} compromises performance by only $2.34\%$; see Fig.~\ref{fig.perf_blk}. We recall that the fully-decentralized controller with structure shown in Fig.~\ref{fig.Gc_M} degrades performance by $3.02\%$; cf.\ Fig.~\ref{fig.sp_ne_coc}. Since the block-sparse controller has more degrees of freedom than the elementwise sparse controller, performance improvement does not come as a surprise. We finally note that the jumps in the number of non-zero elements in Fig.~\ref{fig.perf_blk} are caused by elimination of the entire off-diagonal rows of the feedback gain $K_{r}$ that acts on states different from relative angles.

	\vspace*{-1ex}
\subsection{Comparison of open- and closed-loop systems}
	\label{sec.compare}

We next compare performance of the open-loop system and the closed-loop systems with optimal centralized and fully-decentralized sparse and block-sparse controllers. The structures of these fully-decentralized controllers are shown in Fig.~\ref{fig.Gc_M} and Fig.~\ref{fig.G50}, respectively.

\begin{figure}[]
\centering
\begin{tabular}{rc}
\begin{rotate}{90}{\quad \quad \quad \quad \quad \quad \quad
$\mathrm{Im} \, (\lambda_i)$}\end{rotate}
\!\!\! & \!\!\!
\subfloat
{\includegraphics[width=0.40\textwidth]{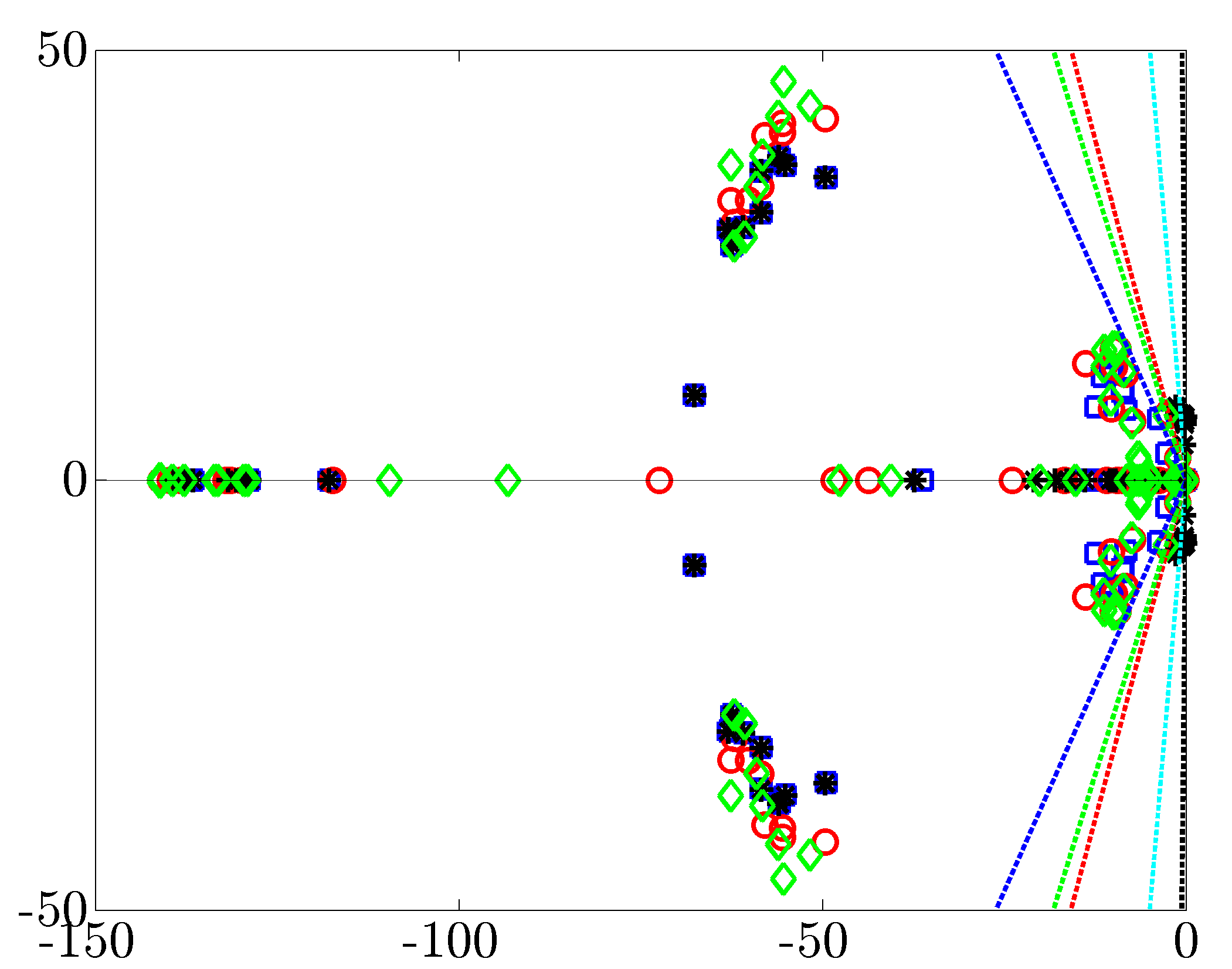}
\label{fig.spectrum_all}}
\\
&
$\mathrm{Re} \, (\lambda_i)$
\\
& (a) Spectra of the open- and close-loop systems
\\[0.25cm]
\begin{rotate}{90}{\quad \quad \quad \quad \quad \quad \quad
$\mathrm{Im} \, (\lambda_i)$}\end{rotate}
\!\!\! & \!\!\!
\subfloat
{\includegraphics[width=0.40\textwidth]{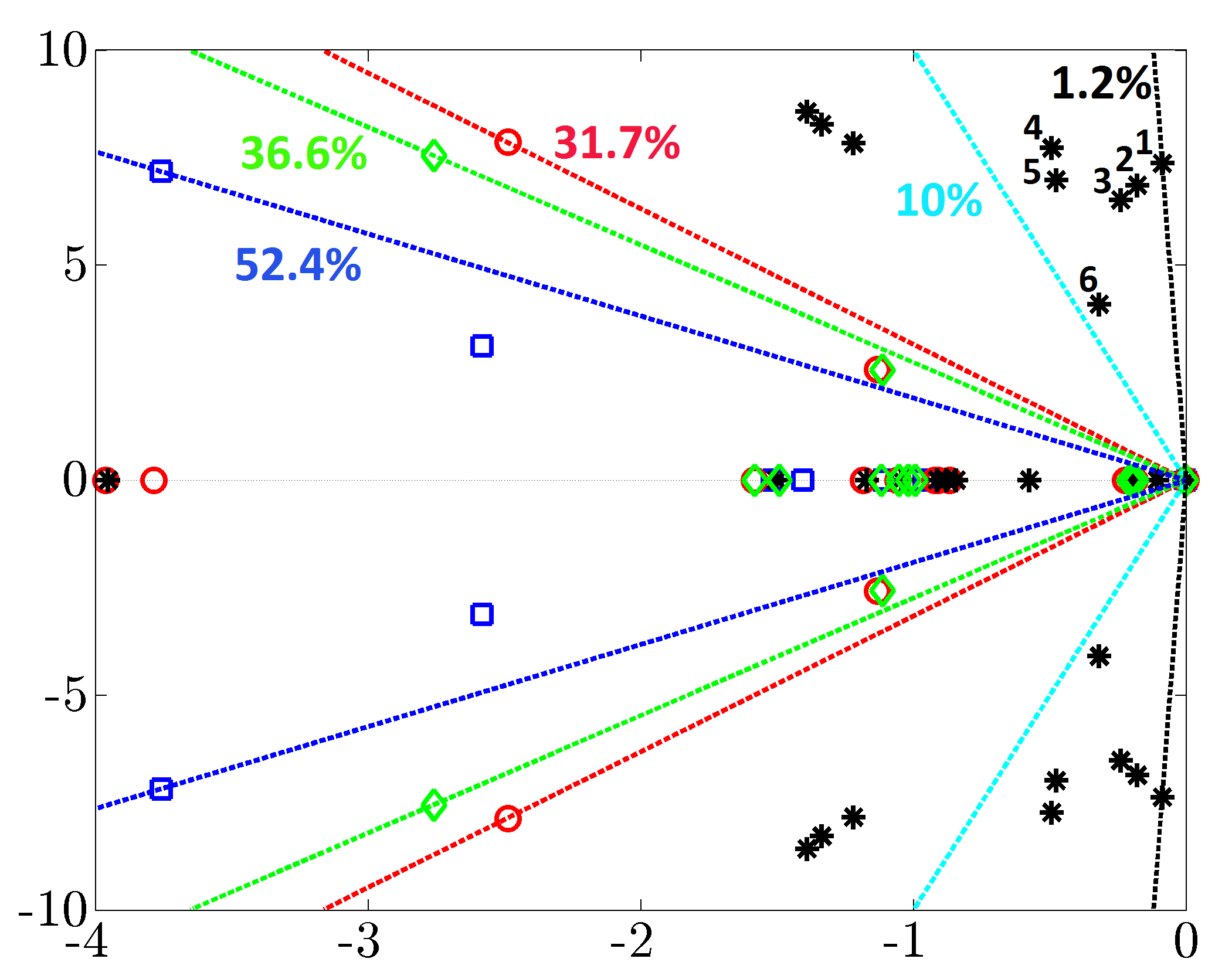}
\label{fig.spectrum_zoom}}
\\
&
$\mathrm{Re} \, (\lambda_i)$
\\
& (b) Zoomed version of (a)
\end{tabular}
\caption{The eigenvalues of the open-loop system and the closed-loop systems with sparse/block-sparse/centralized controllers are represented by $\tc{black}{*}$, $\tc{red}{\circ}$, $\tc{green}{\diamond}$, and $\tc{blue}{\Box}$, respectively. The damping lines indicate lower bounds for damping ratios and they are represented by dashed lines using the same colors as for the respective eigenvalues. The $10\%$ damping line is identified by cyan color. The numbered black asterisks correspond to the six poorly-damped modes given in Table~\ref{Table: poorly-damped modes of New England Grid}.}
\label{fig.spectrum}
\end{figure}

Figure~\ref{fig.spectrum} compares the spectra of the open- and closed-loop systems. As Fig.~\ref{fig.spectrum_all} illustrates, all three controllers (centralized as well as decentralized sparse and block-sparse) move the open-loop spectrum away from the imaginary axis. The dashed lines in Fig.~\ref{fig.spectrum} identify damping lines. Typically, the mode is considered to have sufficient damping if it is located to the left of the $10\%$ cyan damping line. The numbered black asterisks to the right of the $10\%$ damping line in Fig.~\ref{fig.spectrum_zoom} correspond to the six poorly-damped modes of the open-loop system. Other damping lines show that all of our controllers significantly improve the damping of the system by moving the poorly-damped modes deeper into the left-half of the complex plane. This demonstrates that minimization of the variance amplification (i.e., the closed-loop $\mathcal H_{2}$ norm) represents an effective means for improving damping in power systems.

\begin{figure}[]
\centering
\begin{tabular}{c}
\includegraphics[width=0.40\textwidth]{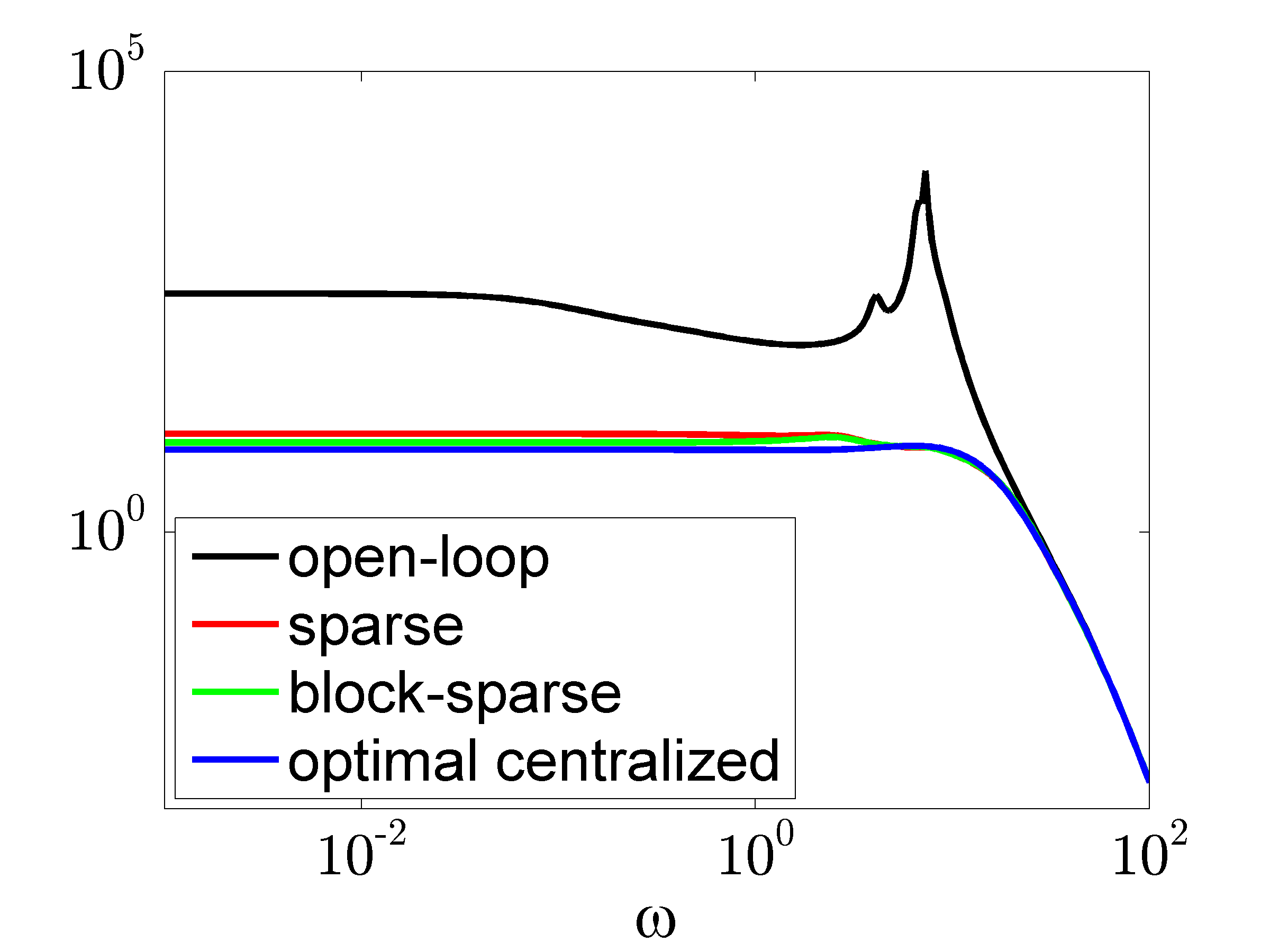}
\end{tabular}
\vspace*{-0.2cm}
\caption{Power spectral density comparison.}
\label{fig.psd}
\end{figure}

Figure~\ref{fig.psd} provides a comparison between the power spectral densities of the four cases. All three controllers successfully suppress the resonant peaks associated with the poorly-damped modes and significantly improve performance. We also note that the fully-decentralized sparse controllers perform almost as well as the optimal centralized controller for high frequencies; for low frequencies, we observe minor discrepancy that accounts for $2-3\%$ of performance degradation in the variance amplification.

\begin{figure}[]
\centering
\begin{tabular}{c}
{eigenvalues of the matrix $Z_1$:}
\\[-0.25cm]
\subfloat[Variance amplification]
{\includegraphics[width=0.3\textwidth]{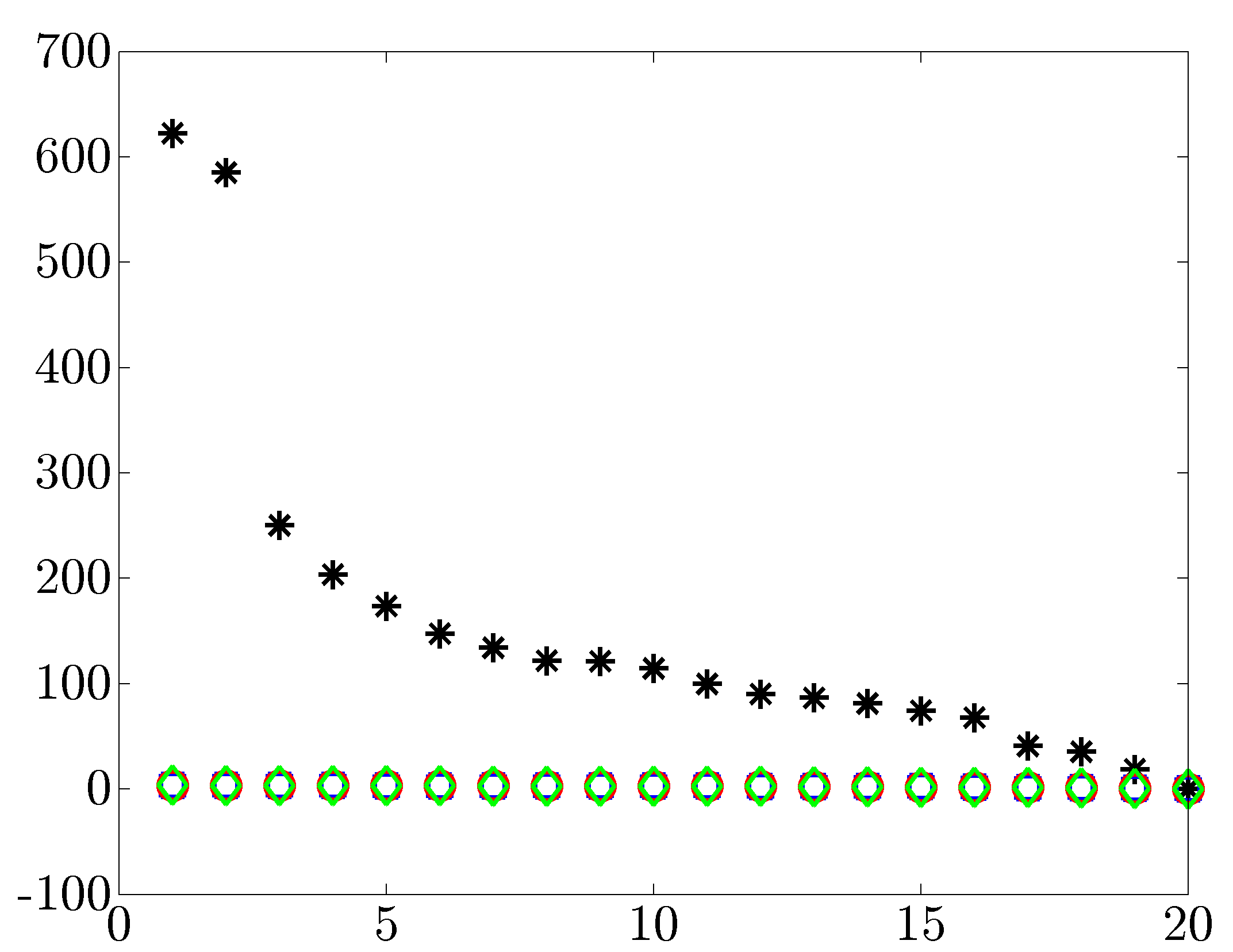}
\label{fig.va_all}}
\\
\subfloat[Zoomed version of (a)]
{\includegraphics[width=0.3\textwidth]{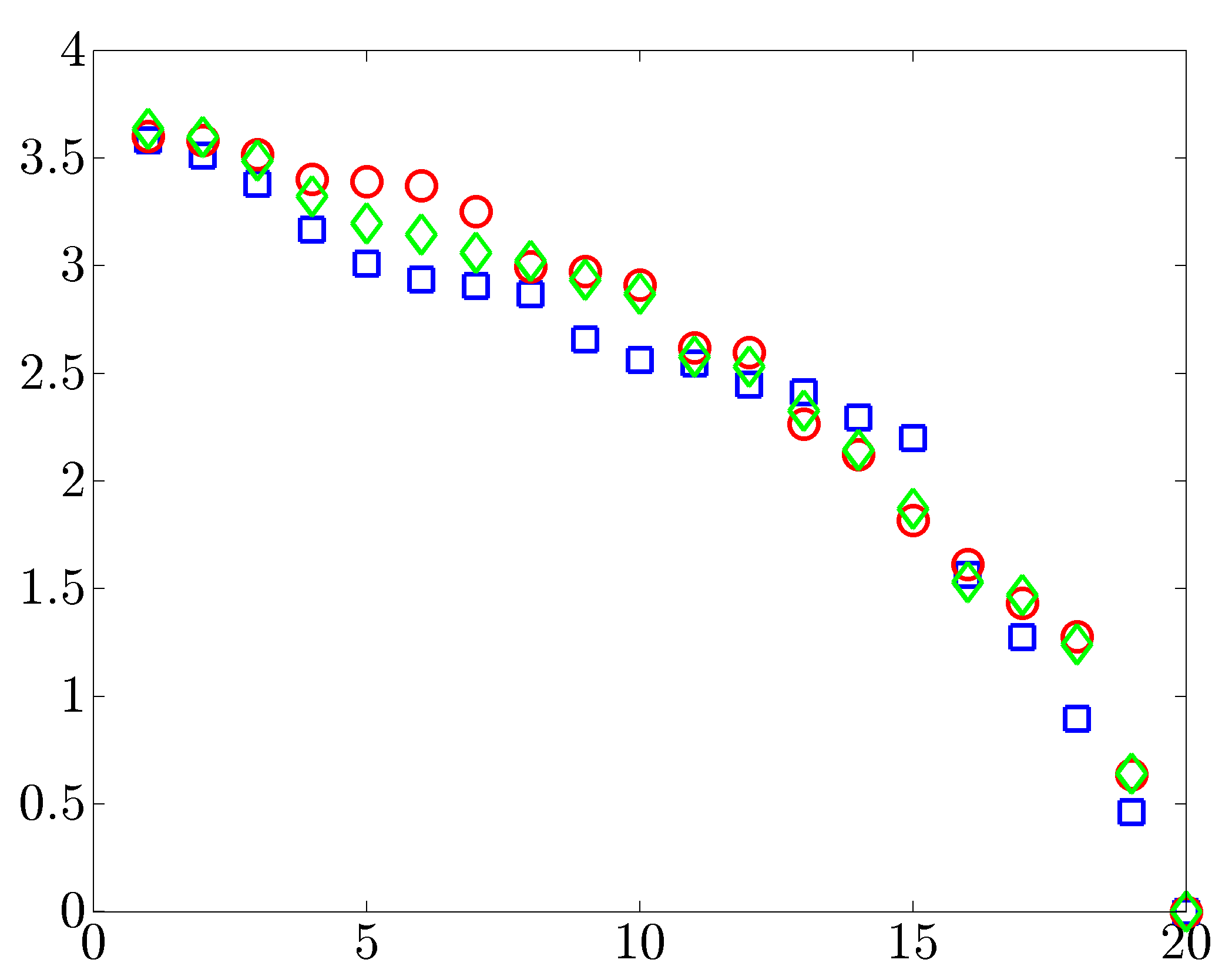}
\label{fig.va_closed}}
\end{tabular}
\caption{Eigenvalues of the output covariance matrix $Z_1$. $*$ represents the open-loop system, \tc{red}{$\circ$}, \tc{green}{$\diamond$} and \tc{blue}{$\Box$} represent the closed-loop systems with sparse, block-sparse, and optimal centralized controllers, respectively.}
\label{fig.va_com}
\end{figure}

Figure~\ref{fig.va_com} displays the eigenvalues of the output covariance matrix $Z_1$ for the four cases mentioned above. Relative to the open-loop system, all three feedback strategies significantly reduce the variance amplification. A closer comparison of the closed-loop systems reveals that the diagonal elements of the output covariance matrix are equalized and balanced by both the optimal centralized and the decentralized controllers; see Fig.~\ref{fig.va_closed}. Similar to the modal observations discussed in~\cite{dorjovchebulTPS14}, the optimal sparse and block-sparse feedback gains not only increase the damping of the eigenvalues associated with the inter-area modes, but also structurally distort these modes by rotating the corresponding eigenvectors.

\begin{figure*}[!t]
\centering
\begin{tabular}{c|cc|c}
& & \rm{Mode 2} & \rm{Mode 6}
\\
\hline
\begin{rotate}{90}{\quad \quad \quad \quad open-loop $\theta$}\end{rotate}
\quad & \quad
\begin{rotate}{90} \quad \quad \quad \quad $\theta_i (t) [ \rm{rad} ]$ \end{rotate}
\!\!\!\!\! & \!\!\!\!\!
\includegraphics[width=0.4\textwidth]{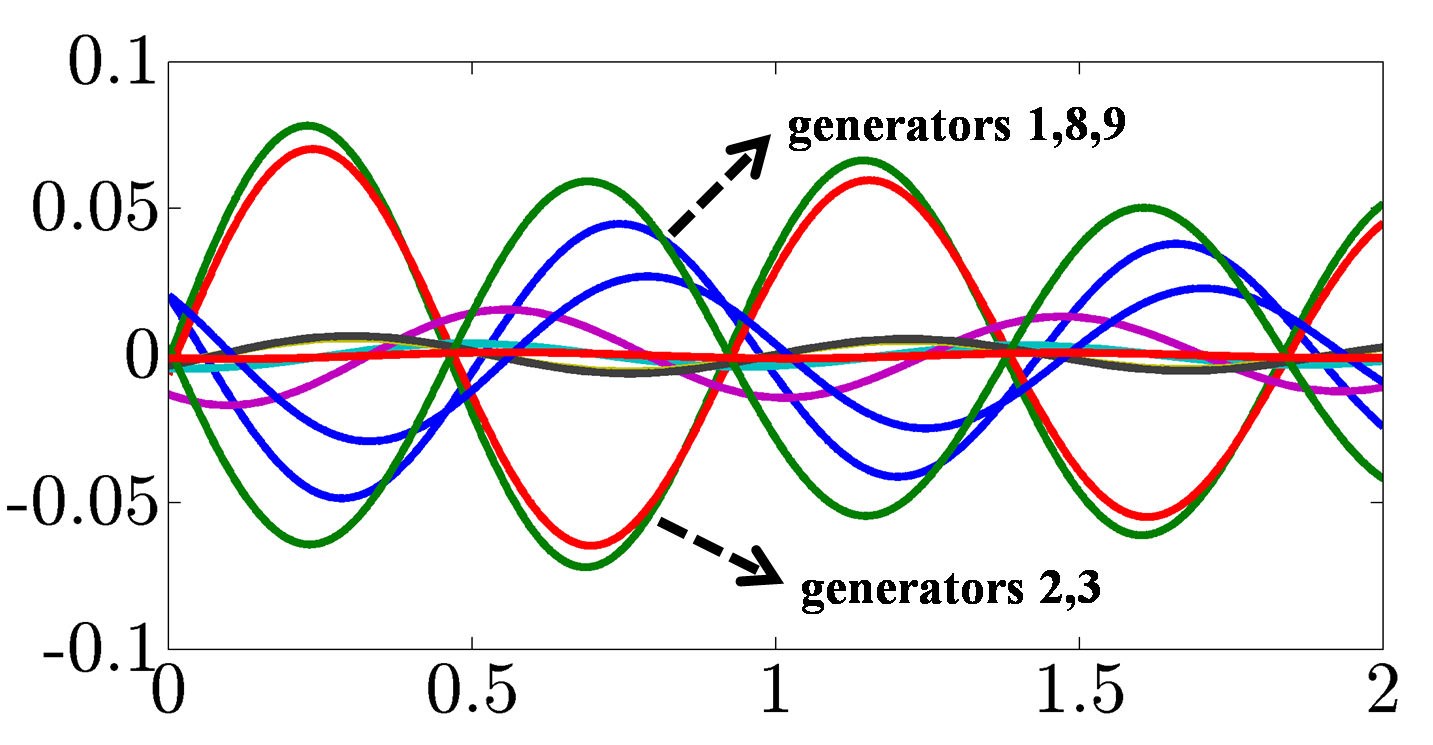}
&
\includegraphics[width=0.4\textwidth]{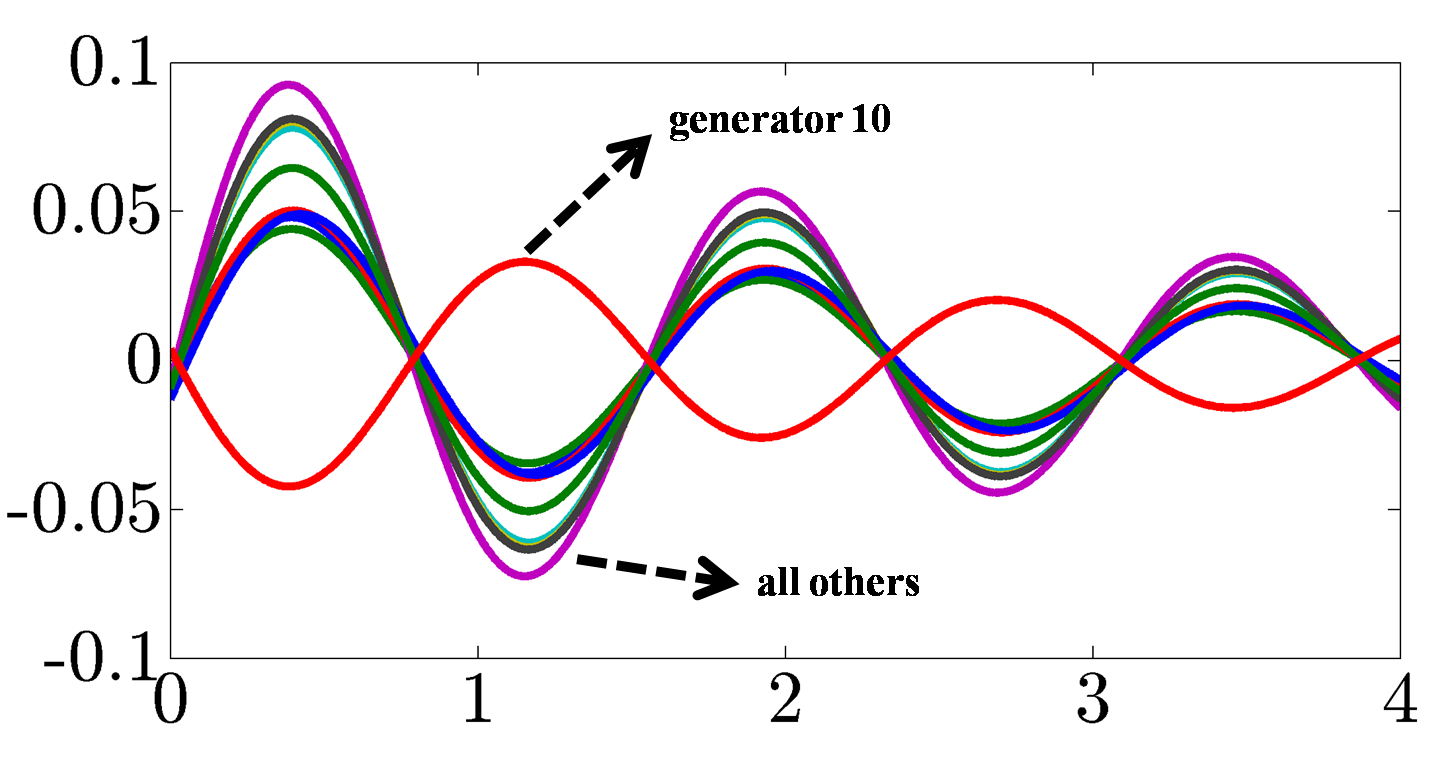}
\\[-0.25cm]
& & $t$ [ \rm{s} ] & $t$ [ \rm{s} ]
\\
\begin{rotate}{90}{\quad \quad \quad closed-loop $\theta$}\end{rotate}
\quad & \quad
\begin{rotate}{90} \quad \quad \quad \quad $\theta_i (t) [ \rm{rad} ]$ \end{rotate}
\!\!\!\!\! & \!\!\!\!\!
\includegraphics[width=0.4\textwidth]{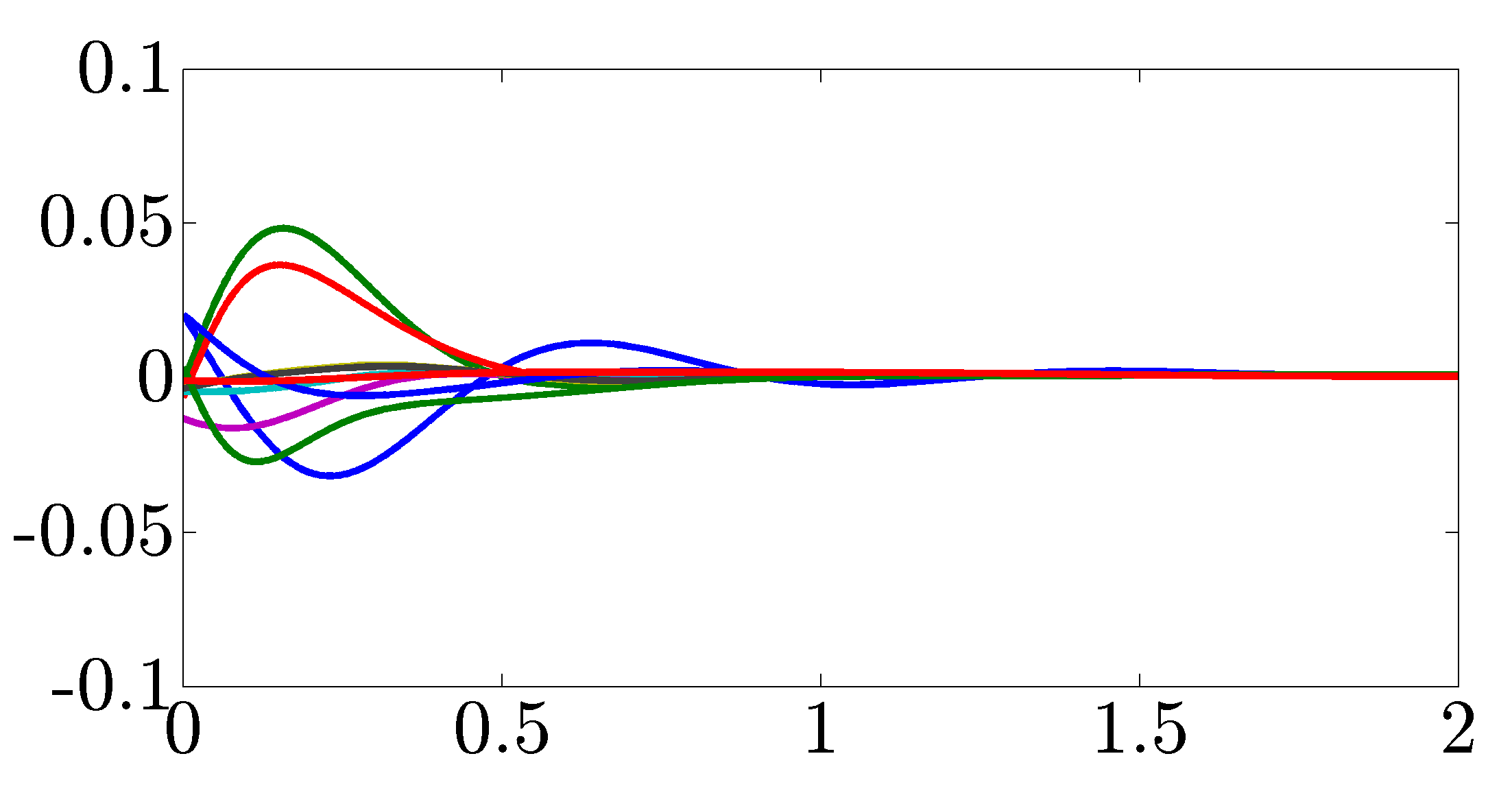}
&
\includegraphics[width=0.4\textwidth]{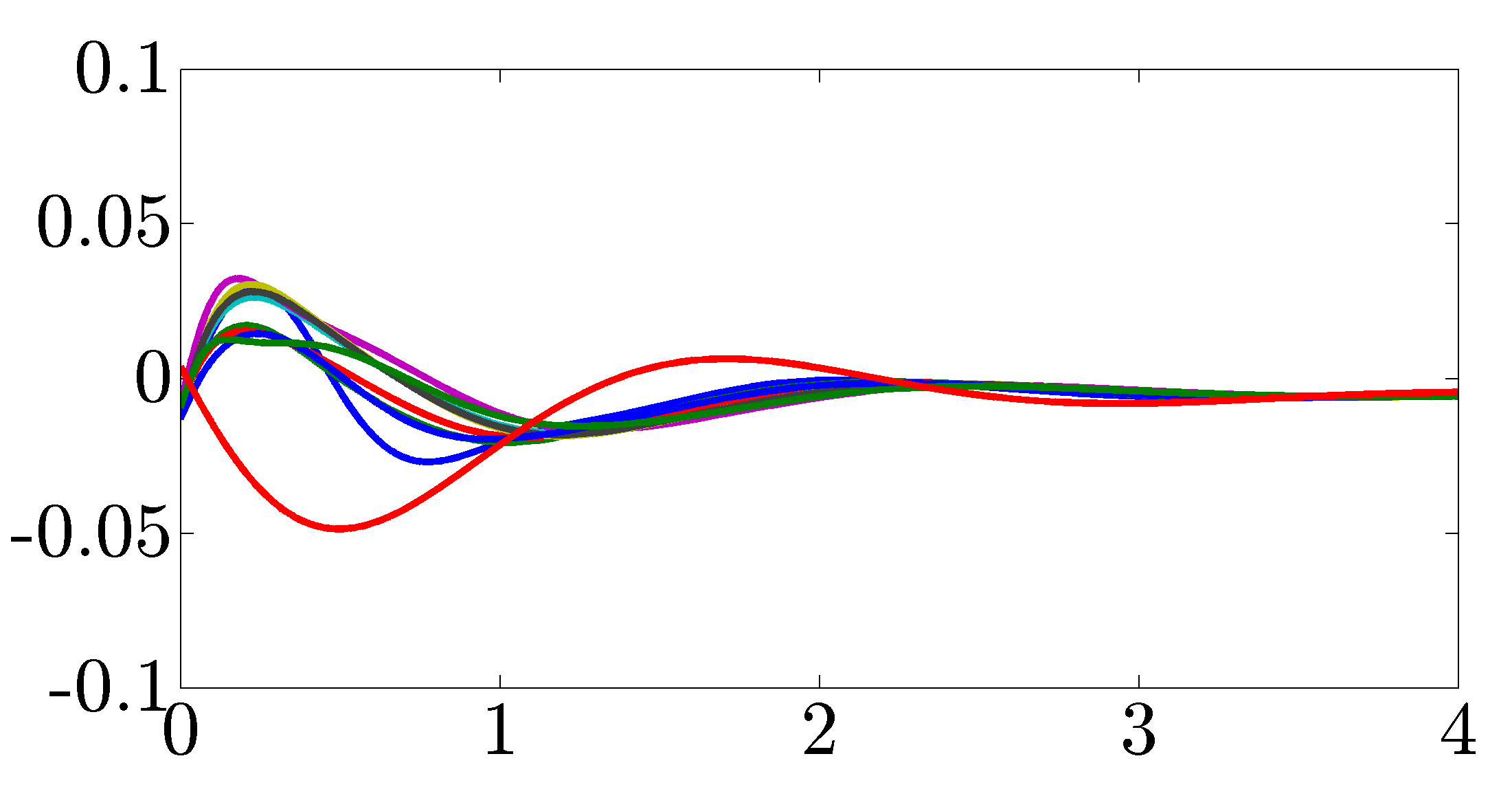}
\\[-0.25cm]
& & $t$ [ \rm{s} ] & $t$ [ \rm{s} ]
\\
\hline
\begin{rotate}{90}{\quad \quad \quad \quad open-loop $\dot{\theta}$}\end{rotate}
\quad & \quad
\begin{rotate}{90} \quad \quad \quad $\dot{\theta}_i (t) [ \rm{rad/s} ]$ \end{rotate}
\!\!\!\!\! & \!\!\!\!\!
\includegraphics[width=0.4\textwidth]{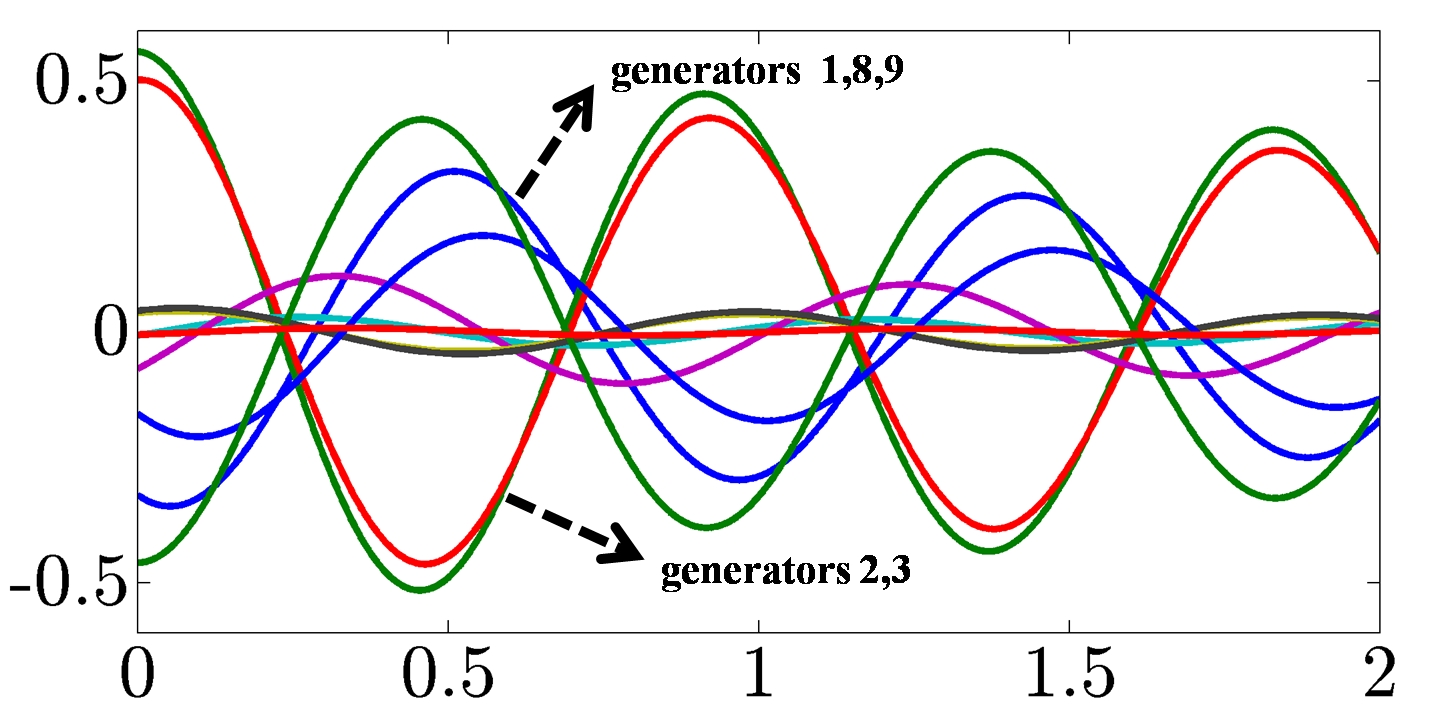}
&
\includegraphics[width=0.4\textwidth]{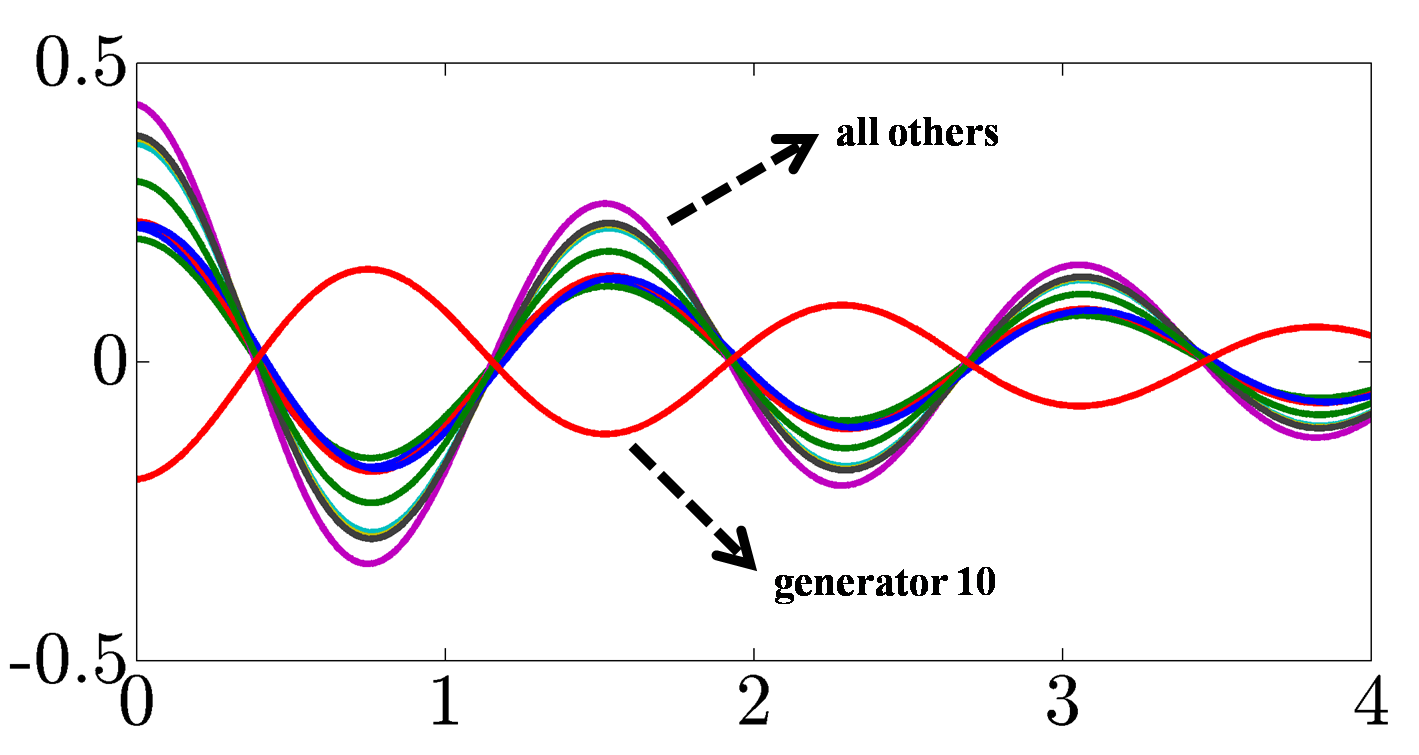}
\\[-0.25cm]
& & $t$ [ \rm{s} ] & $t$ [ \rm{s} ]
\\
\begin{rotate}{90}{\quad \quad \quad closed-loop $\dot{\theta}$}\end{rotate}
\quad & \quad
\begin{rotate}{90} \quad \quad \quad $\dot{\theta}_i (t) [ \rm{rad/s} ]$ \end{rotate}
\!\!\!\!\! & \!\!\!\!\!
\includegraphics[width=0.4\textwidth]{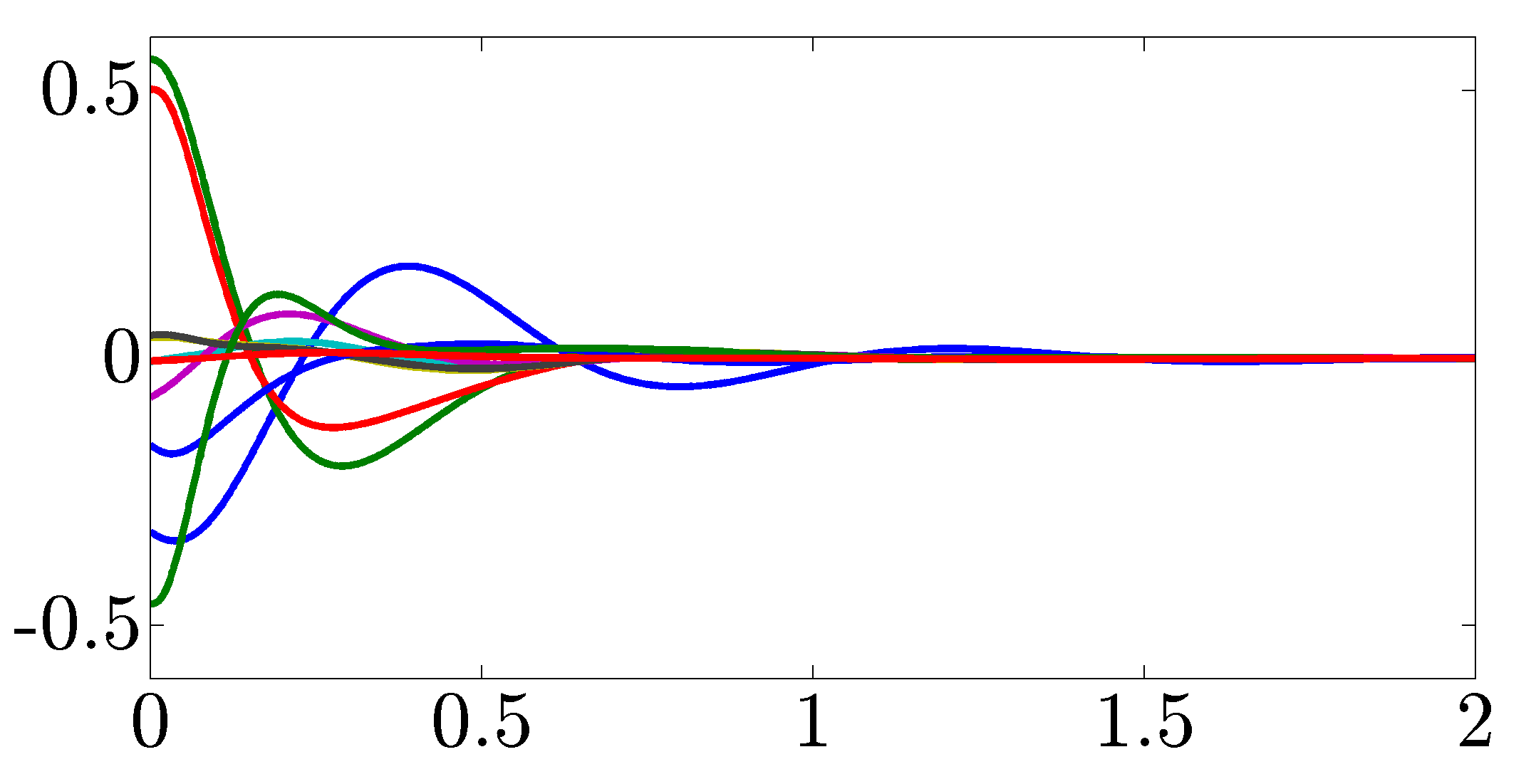}
&
\includegraphics[width=0.4\textwidth]{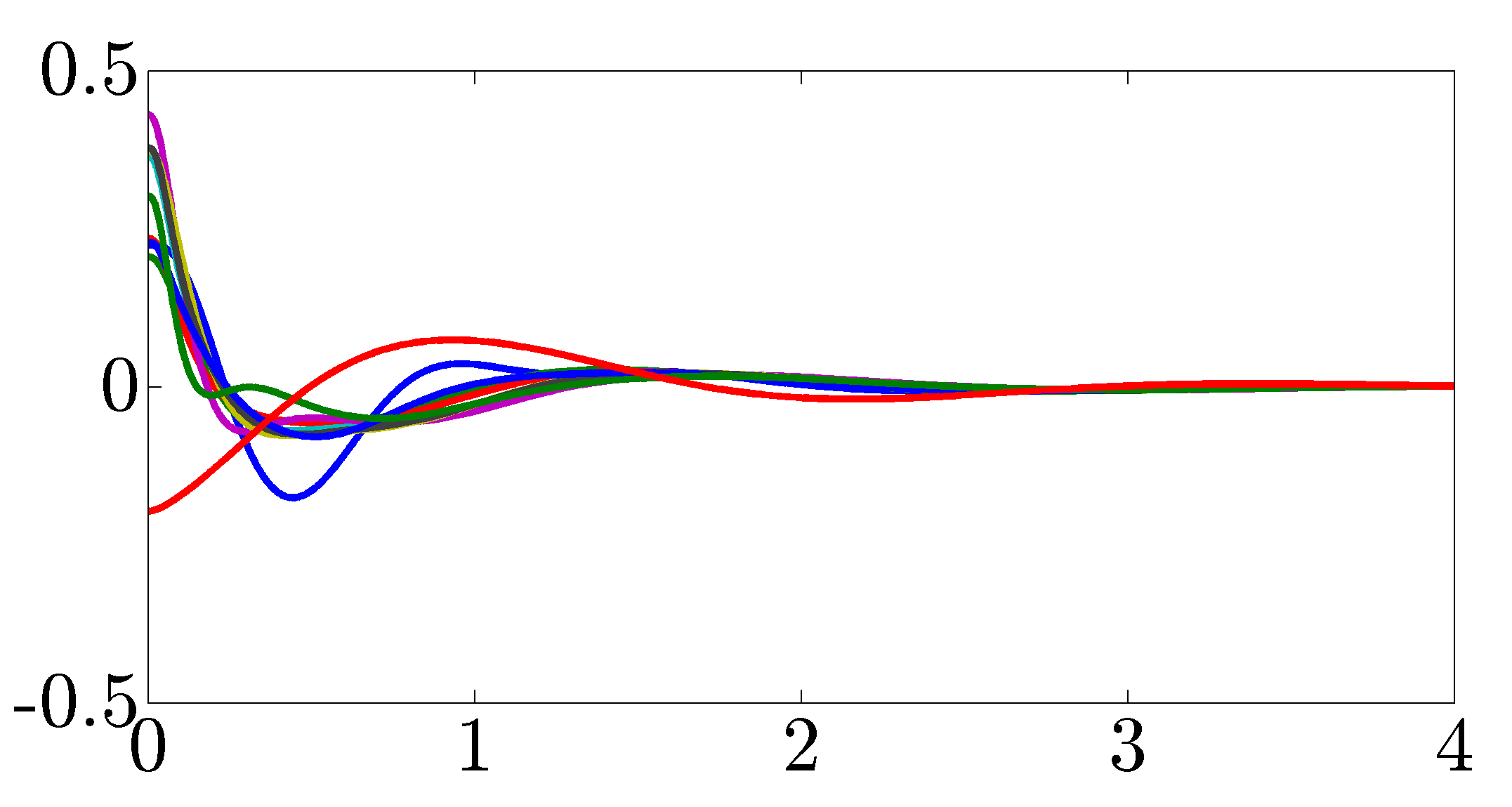}
\\[-0.25cm]
& & $t$ [ \rm{s} ] & $t$ [ \rm{s} ]
\end{tabular}
\caption{Time-domain simulations of the linearized model of the IEEE 39 New England power grid. The rotor angles and frequencies of all generators are shown. The closed-loop results are obtained using the fully-decentralized block-sparse controller. The initial conditions are given by the eigenvectors of the poorly-damped inter-area modes $2$ (left) and $6$ (right) from Table~\ref{Table: poorly-damped modes of New England Grid}.}
\label{fig.simulation}
\end{figure*}

We use time-domain simulations of the linearized model to verify performance of decentralized block-sparse controller. Figure~\ref{fig.simulation} shows the trajectories of rotor angles and frequencies for the open- and closed-loop systems for two sets of initial conditions. These are determined by the eigenvectors of open-loop inter-area modes $2$ and $6$ in Table~\ref{Table: poorly-damped modes of New England Grid}. Clearly, the decentralized block-sparse controller significantly improves performance by suppressing the inter-area oscillations between groups of generators. Furthermore, relative to the open-loop system, the transient response of the closed-loop system features shorter settling time and smaller maximum overshoot. 

\begin{figure}
\centering
\begin{tabular}{ccc}
\begin{rotate}{90}{\scriptsize \quad $\sharp$ of operating points}\end{rotate}
&
\subfloat[Open-loop system]
{\includegraphics[width=0.20\textwidth]{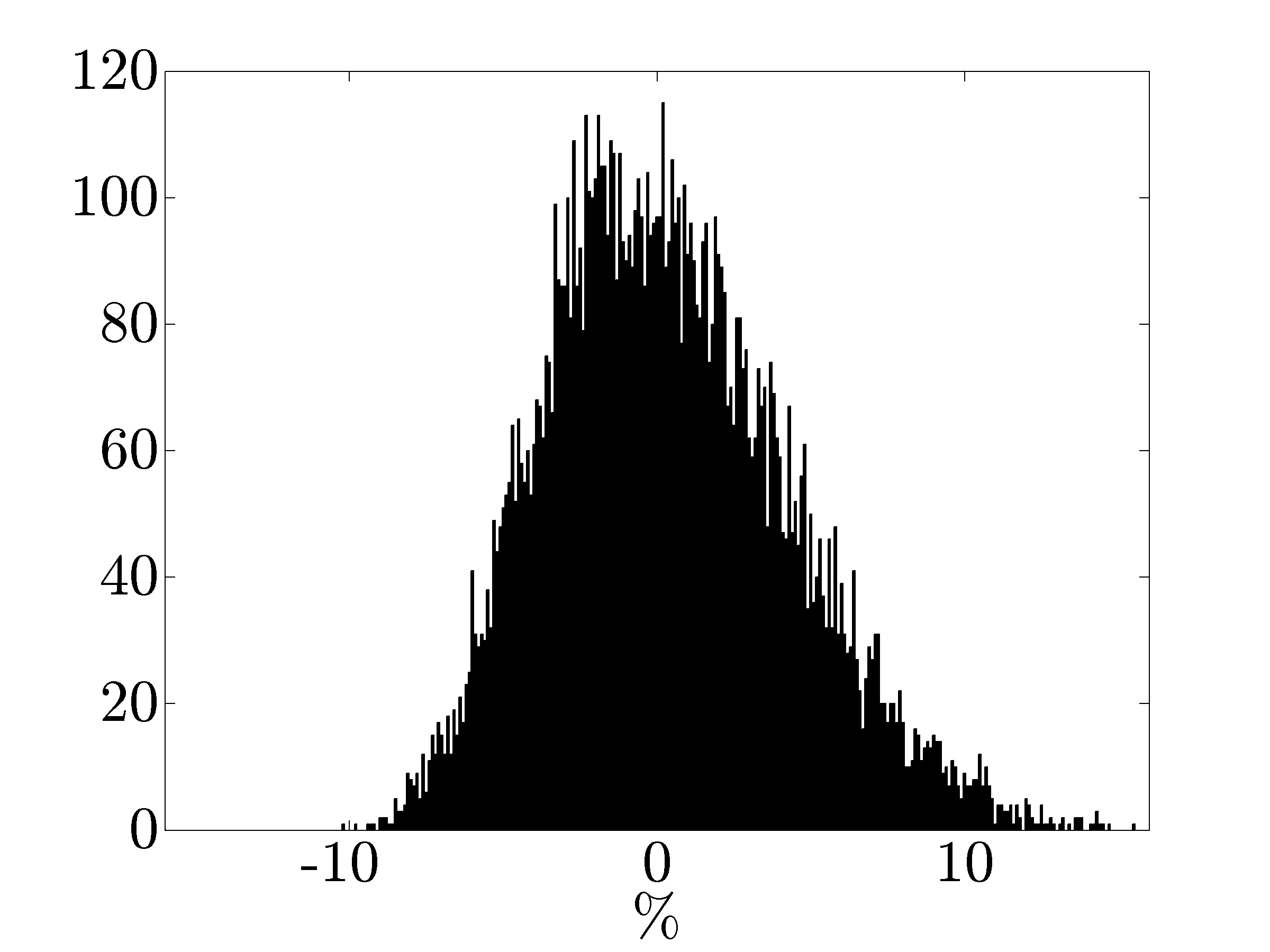}
\label{fig.hist_open}}
&
\subfloat[Centralized controller]
{\includegraphics[width=0.20\textwidth]{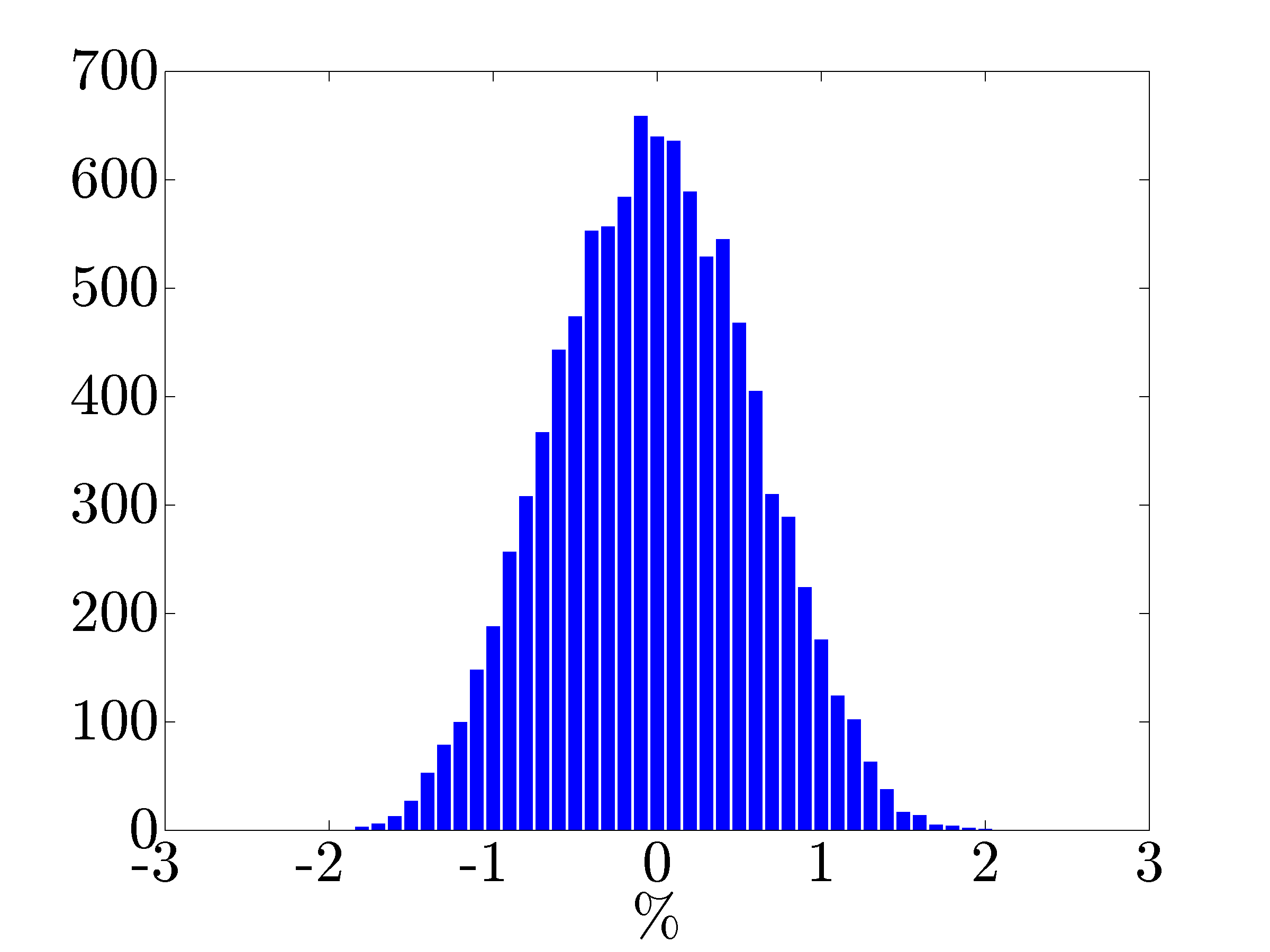}
\label{fig.hist_central}}
\\
\begin{rotate}{90}{\scriptsize \quad $\sharp$ of operating points}\end{rotate}
&
\subfloat[Sparse controller]
{\includegraphics[width=0.20\textwidth]{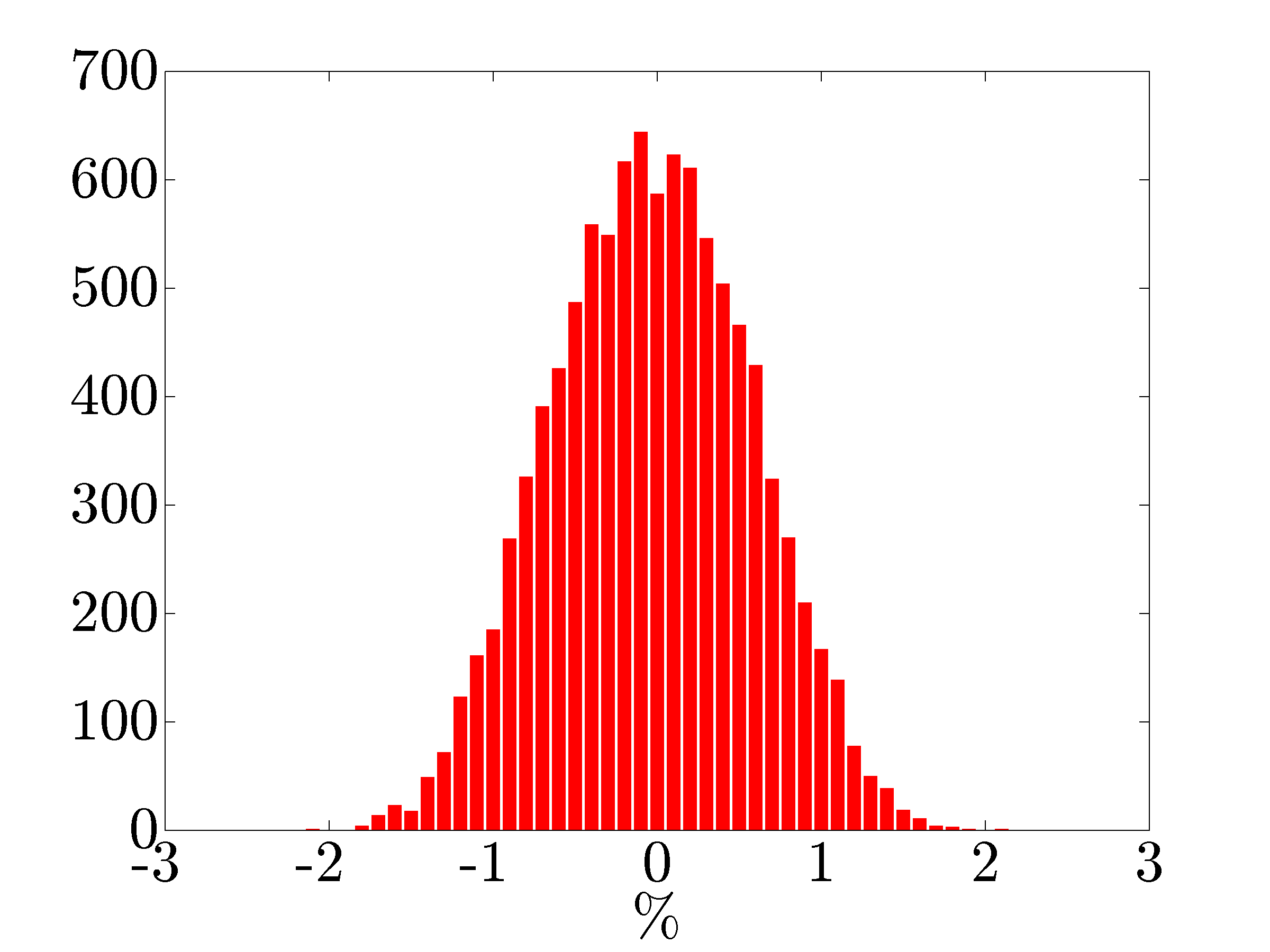}
\label{fig.hist}}
&
\subfloat[Block-sparse controller]
{\includegraphics[width=0.20\textwidth]{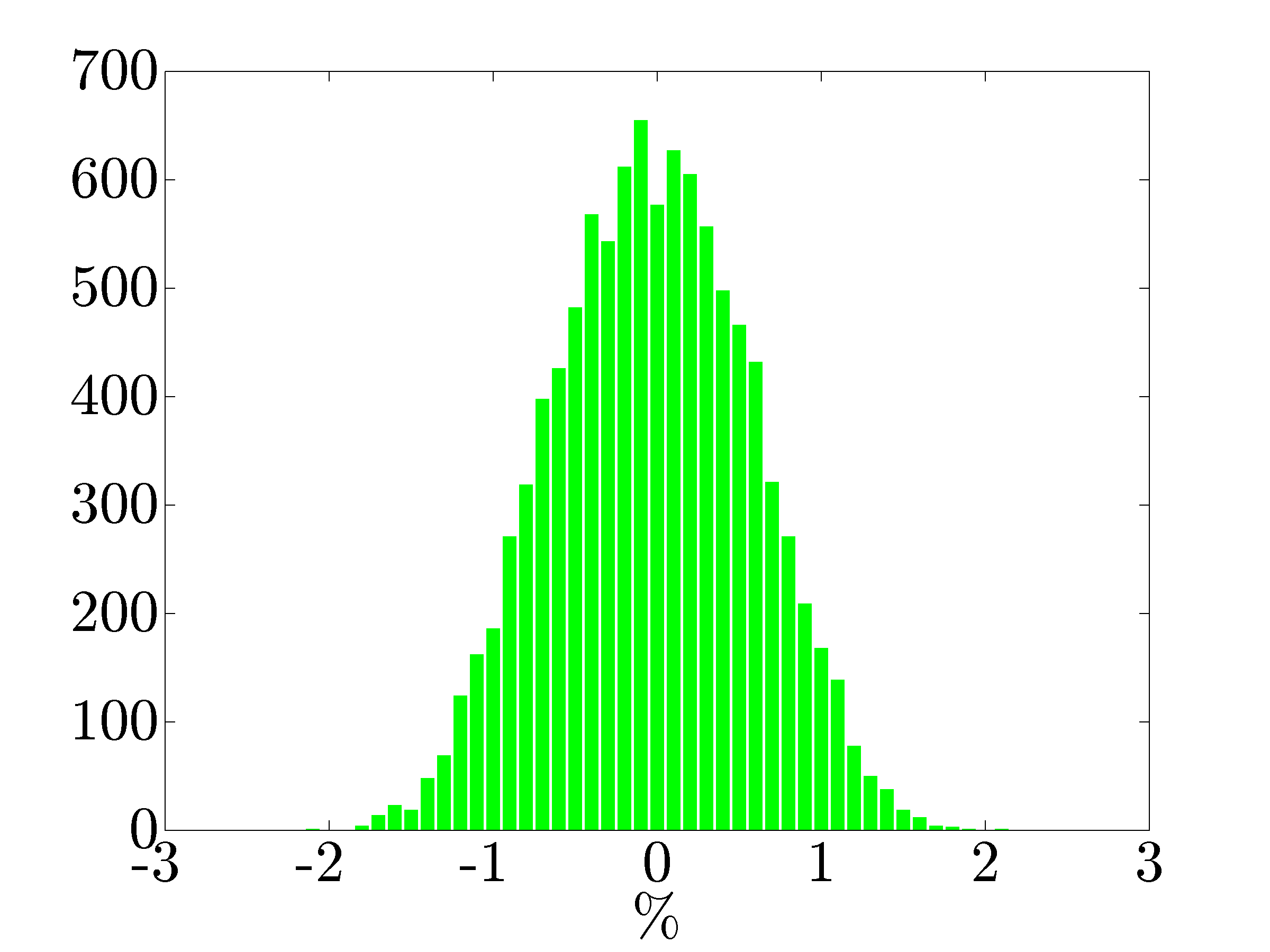}
\label{fig.hist_blk}}
\end{tabular}
\caption{Performance histograms of open- and closed-loop linearized systems (with nominal controllers) for $10,000$ uniformly distributed operating points.}
\label{fig.H2norm}
\end{figure}

\begin{figure}
\centering
\begin{tabular}{cc}
\begin{rotate}{90}{\footnotesize \quad \quad \quad \quad \quad \quad degree}\end{rotate}
&
{\includegraphics[width=0.4\textwidth]{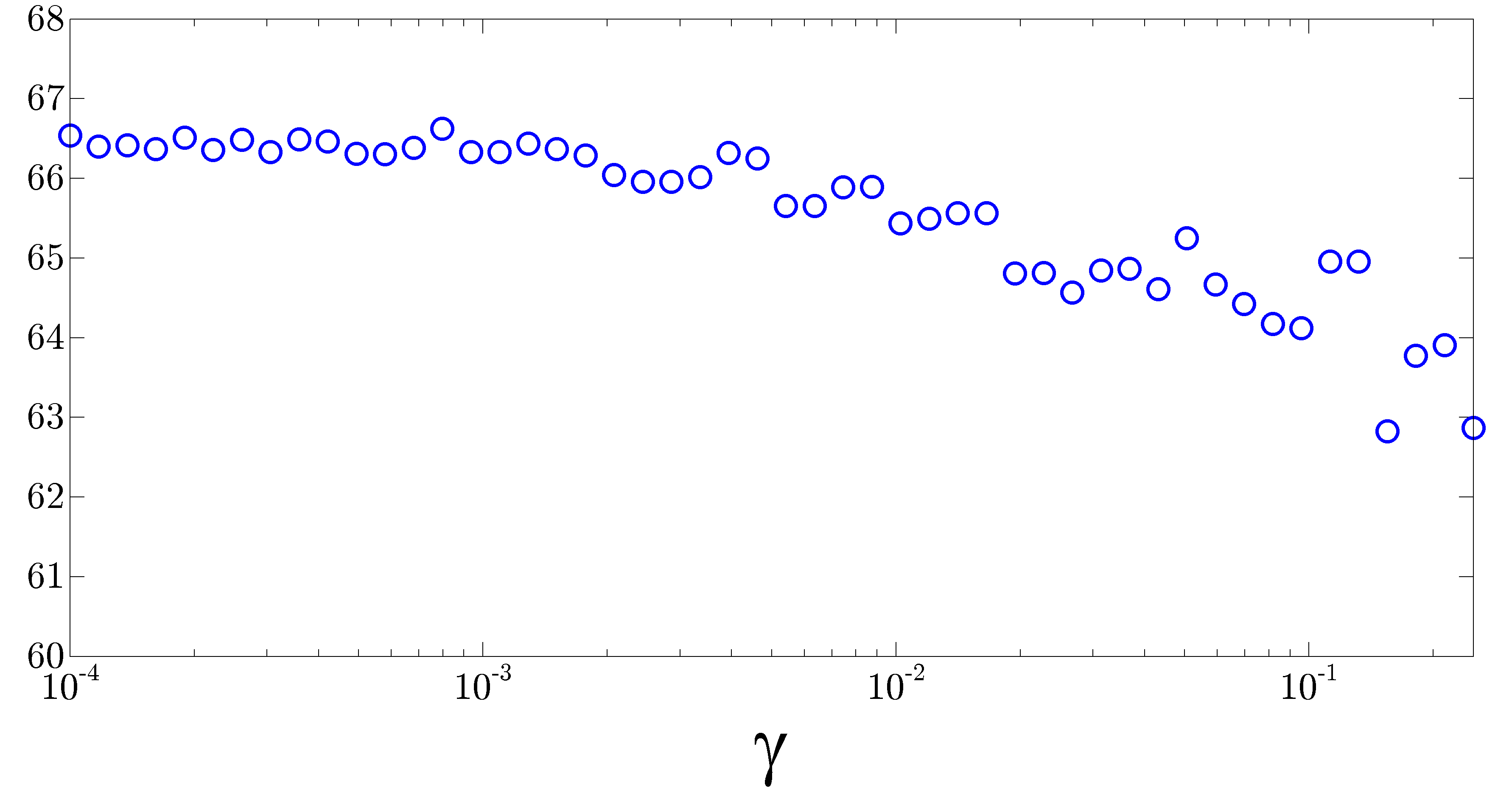}
\label{fig.robustness1}}
\\
& {(a) element-wise sparse controller}
\\[0.25cm]
\begin{rotate}{90}{\footnotesize \quad \quad \quad \quad \quad \quad degree}\end{rotate}
&
{\includegraphics[width=0.4\textwidth]{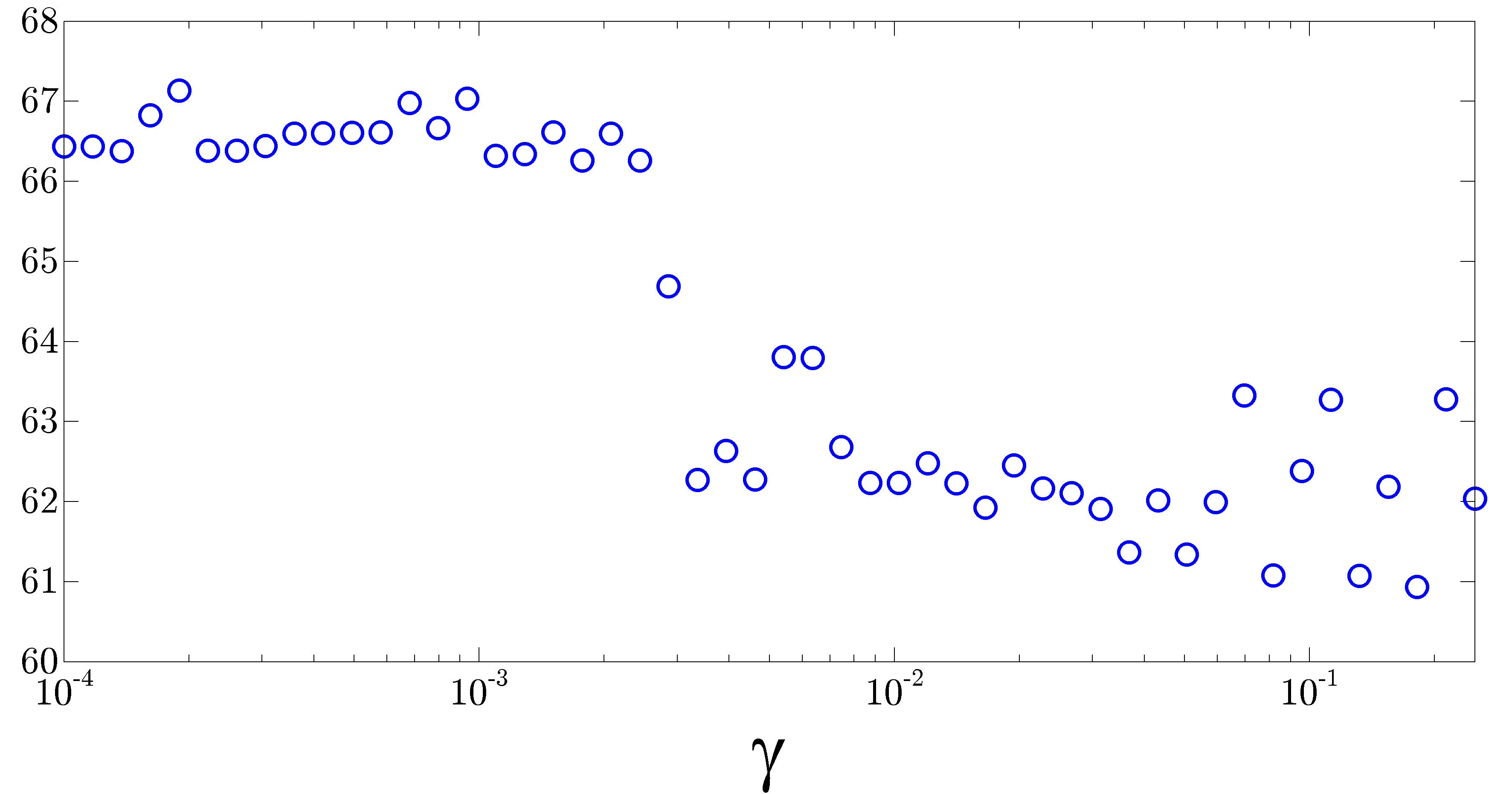}
\label{fig.robustness2}}
\\
& {(b) block-sparse controller}
\end{tabular}
\caption{Multivariable phase margins as a function of $\gamma$.}
\label{fig.robustness}
\end{figure}

\subsection{Robustness analysis}

We close this section by examining robustness to the operating point changes of both open- and closed-loop systems. Random load perturbations are used to modify the operating point of the nonlinear system. The loads, that are used for the analysis and control synthesis, are altered via uniformly distributed perturbations that are within $\pm 20\%$ of the nominal loads. The performance of the {\em nominal\/} centralized and decentralized controllers on the {\em perturbed linearized model\/} is evaluated by examining the closed-loop $\mathcal H_{2}$ norm.

Figure~\ref{fig.H2norm} shows the distribution of performance change for $10,000$ operating points around the original equilibria. We observe bell-shaped distributions with symmetric and narrow spread around the nominal performance. In spite of significant changes in the operating points, both centralized and fully-decentralized controllers are within $2 \%$ of the nominal performance. In contrast, same perturbations can degrade performance of the open-loop system by as much as $15 \%$. Thus, our decentralized controllers also reduce the sensitivity and improve the robustness with respect to setpoint changes.

To account for delays in communication channels, asynchronous measurements, and fast unmodeled dynamics, we utilize multivariable phase margin to quantify the robustness of our sparse optimal controllers. In Fig.~\ref{fig.robustness}, we investigate how the phase margins of the closed-loop systems change with the sparsity-promoting parameter $\gamma$. As our emphasis on sparsity increases, multivariable phase margins degrade gracefully and stay close to a desirable phase margin of $60^{\circ}$.

Our approach thus provides a systematic way for designing optimal sparse controllers with favorable robustness margins and performance guarantees even in a fully-decentralized case.

\section{Concluding remarks}
\label{Section:Conclusion}

We have analyzed inter-area oscillations in power systems by studying their power spectral densities and output covariances. Our analysis of the open-loop system identifies poorly-damped modes that cause inter-area oscillations. We have also designed sparse and block-sparse feedback controllers that use relative angle measurements to achieve a balance between system performance and controller architecture. By placing increasing weight on the sparsity-promoting term we obtain fully-decentralized feedback gains. Performance comparisons of open- and closed-loop systems allowed us to understand the effect of the control design approach both in terms of system performance and with regards to the resulting control architecture. For the IEEE 39 New England model we have successfully tested our analysis and control design algorithms. We have also provided a systematic method for optimal retuning of fully-decentralized excitation controllers that achieves comparable performance to the optimal centralized controller.

\vspace*{-2ex}
\renewcommand{\baselinestretch}{1}

\end{document}